\documentclass[11pt,a4paper,twoside]{article}
\usepackage{latexsym,amsfonts,amsmath, amsthm, amssymb}
\usepackage{fullpage}
\usepackage{xcolor,bm}
\usepackage[utf8x]{inputenc}

\newtheorem{theorem}{Theorem}
\newtheorem{proposition}[theorem]{Proposition}
\newtheorem{lemma}[theorem]{Lemma}
\newtheorem{corollary}[theorem]{Corollary}

\theoremstyle{definition}

\newtheorem{remark}[theorem]{Remark}

\newcommand{\ls}{\leqslant}
\newcommand{\gs}{\geqslant}

\begin{document}

\pagecolor{white}

\title{\bf Regular ellipsoids and a Blaschke-Santaló-type inequality for projections of non-symmetric convex bodies}

\medskip

\author{Beatrice-Helen Vritsiou}

\date{}

\maketitle

\begin{abstract}
\small 
It is shown that every not-necessarily symmetric convex body $K$ in ${\mathbb R}^n$
has an affine image $\tilde{K}$ of $K$ 
such that the covering numbers of $\tilde{K}$ by growing dilates of the unit Euclidean ball,
as well as those of the unit Euclidean ball by growing dilates of $\tilde{K}$, decrease in a regular way.
This extends to the non-symmetric case a famous theorem by Pisier, albeit with worse estimates 
on the rate of decrease of the covering numbers.
The affine image $\tilde{K}$ can be chosen to have either barycentre or Santaló point at the origin.

In the proof we use Pisier's theorem as a black box, as well as a suggested approach by Klartag and V. Milman.
A key new ingredient is Blaschke-Santal\'{o}-type inequalities for projections of a body $\tilde{K}$ 
with Santal\'{o} point at the origin, which could be of independent interest. Unlike the application to covering, these (as well as the analogous inequalities for centred convex bodies that were already considered by Klartag and Milman \cite{Klartag-VMilman-preprint}) can be shown to be optimal up to absolute constants.

We also present an application to results around the mean norm of isotropic (not-necessarily symmetric) convex bodies.
\end{abstract}

\section{Introduction}\label{sec:intro}

Given convex bodies $K,L$ in ${\mathbb R}^n$, the covering number $N(K,L)$ of $K$ by $L$
is the least number of translates of $L$ whose union covers $K$. 
%Clearly this is finite if $K$ is bounded and $L$ has nonempty interior; in this note we will only be concerned
%with convex bodies $K, L$, that is, convex and compact sets that have nonempty interior
%(in the presence of convexity, the latter property is equivalent to $K$ having positive Lebesgue volume 
%which we will denote by $|K|$). 
%
%%Given sets $K,L$ in ${\mathbb R}^n$, the covering number $N(K,L)$ of $K$ by $L$
%%is the least number of translates of $L$ whose union covers $K$. Clearly this is finite
%%if $K$ is bounded and $L$ has nonempty interior; in this note we will only be concerned
%%with convex bodies $K, L$, that is, convex and compact sets that have nonempty interior
%%(in the presence of convexity, the latter property is equivalent to $K$ having positive Lebesgue volume 
%%which we will denote by $|K|$). 
%
Let us initially assume that $K$ is origin-symmetric (that is, $K=-K$), and that $L$ is the unit Euclidean ball $B_2^n$
in ${\mathbb R}^n$.  
A celebrated theorem of V. Milman \cite{VMilman-1986} states that 
any such 
%convex body 
$K$ has a linear image $\tilde{K} = T(K)$, $T\in {\rm GL}(n)$, such that 
\begin{equation} \label{Milman-Mposition}
\max\{N(\tilde{K}, B_2^n),N(\tilde{K}^\circ, B_2^n), N(B_2^n, \tilde{K}), N(B_2^n, \tilde{K}^\circ)\} \ls \exp(Cn) \end{equation}
for some constant $C>0$ independent of $n$ or $K$. Here $L^\circ$ denotes the polar set of a convex body $L$ which contains $0$ in its interior:
%containing the origin in its interior: 
%$L^\circ := \{y\in {\mathbb R}^n: \langle x,y\rangle\ls 1\ \hbox{for all}\ x\in L\}$.
\begin{equation*} 
L^\circ := \{y\in {\mathbb R}^n: \langle x,y\rangle\ls 1\ \hbox{for all}\ x\in L\}.
\end{equation*}
%The estimate is optimal up to the value of the absolute constant $C$, as can be seen by simple volumetric considerations:
%$N(K,L) \gs |K|/|L|$, where $|A|\equiv |A|_n$ denotes the Lebesgue volume of a set $A\subseteq {\mathbb R}^n$. 
%On the other hand, it does not offer information on how the covering numbers may decrease 
%if we consider larger and larger dilates of the covering body.

Pisier \cite{Pisier-1989} refined this theorem by acquiring additional information
on how the covering numbers may decrease if we consider larger and larger dilates of the covering bodies.

\begin{theorem}
{\rm (Pisier, 1989)} \label{thm:Pisier-regular-ellipsoids}
For every $\alpha \in (0,2)$ there exists a constant $C_\alpha \gs 1$ such that
the following holds: for every symmetric convex body $K$ in ${\mathbb R}^n$ we can find
an ellipsoid ${\cal E} = {\cal E}(K,\alpha)$ in ${\mathbb R}^n$ with the property that, for every $t\gs C_\alpha^{1/\alpha}$, 
\begin{equation}\label{Pisier-Mposition1}
\max\left\{N(K, t{\cal E}), N(K^\circ, t{\cal E}^\circ), N({\cal E}, tK), N({\cal E}^\circ, tK^\circ)\right\} 
\ls \exp\bigl(C_\alpha\, n/t^\alpha\bigr).
\end{equation}
%for every $t\geq 1$. 
Equivalently, 
%for every symmetric convex body $K$ in ${\mathbb R}^n$ 
we can find $T\in {\rm GL}(n)$ such that, setting $\tilde{K} = T(K)$, we will have
\begin{equation}\label{Pisier-Mposition2}
\max\{N(\tilde{K}, t B_2^n),N(\tilde{K}^\circ, t B_2^n), N(B_2^n, t\tilde{K}), N(B_2^n, t\tilde{K}^\circ)\}
\ls \exp\bigl(C_\alpha\, n/t^\alpha\bigr)
\end{equation}
for every $t\gs C_\alpha^{1/\alpha}$.

The constants $C_\alpha$ satisfy
%$C_\alpha = O\Bigl(\left(\frac{\alpha}{2-\alpha}\right)^{\alpha/2}\Bigr)$ 
$C_\alpha = O\Bigl(\left(2-\alpha\right)^{-\alpha/2}\Bigr)$ 
as $\alpha \to 2^-$.
\end{theorem}

We say that a body $\tilde{K}$ that satisfies \eqref{Milman-Mposition} is in $M$-position,
or that $\tilde{K} = T(K)$ is the $M$-position of $K$.
Similarly, we say that a body $\tilde{K}$ that satisfies \eqref{Pisier-Mposition2}
for some $\alpha\in (0,2)$ is in $\alpha$-regular $M$-position,
while the ellipsoid ${\cal E}(K,\alpha)$ appearing in \eqref{Pisier-Mposition1} 
is called an $\alpha$-regular $M$-ellipsoid of $K$.

\medskip

In the case of Milman's theorem, it was observed shortly thereafter that, with the use of other central results in Convex Geometry, it can be extended to non-necessarily-symmetric convex bodies $K$: the most usual choices for positioning $K$ are to assume its barycentre or its Santaló point is at the origin, and then some linear transformation brings it to $M$-position (we briefly go over the details for this in Section \ref{sec:prelim}).

%Note that \eqref{Milman-Mposition} is optimal up to the value of the absolute constant $C$.
%Indeed, by the obvious volumetric bound following from the definition of covering numbers,
%\begin{equation*}
%\frac{|B_2^n|^2}{|K||K^\circ|} = \frac{|B_2^n|^2}{|\tilde{K}||\tilde{K}^\circ|} \ls N(B_2^n, \tilde{K}) \cdot N(B_2^n, \tilde{K}^\circ)
%\ls \bigl(\max\bigr\{N(B_2^n, \tilde{K}), N(B_2^n, \tilde{K}^\circ)\bigr\}\bigr)^2,
%\end{equation*}
%and we know that $|K||K^\circ|$, which is a linear invariant of $K$, can differ from $|B_2^n|^2$
%(here $|A|$ denotes the $n$-dimensional Lebesgue volume of a Borel set $A\subseteq {\mathbb R}^n$).

In the case of Pisier's theorem though, the reduction of the non-symmetric case to the already known symmetric case was much more delicate. This was the main motivation of this paper:

\begin{theorem}\label{main-result}
For every $\beta \in \bigl(0,\frac{2}{5}\bigr)$ there exists a constant $D_\beta \gs 1$ such that
the following holds: for every convex body $K$ in ${\mathbb R}^n$
we can find an affine image $\tilde{K}$ of $K$ (namely a linear transformation of a translate of $K$) which satisfies
\begin{equation}\label{general-Mposition}
\max\{N(\tilde{K}, t B_2^n),N(\tilde{K}^\circ, t B_2^n), N(B_2^n, t\tilde{K}), N(B_2^n, t\tilde{K}^\circ)\}
\ls \exp\bigl(D_\beta\, n/t^\beta\bigr)
\end{equation}
for every $t\gs D_\beta^{1/\beta}$.

Moreover, if $K$ has barycentre or Santaló point at the origin, then
only a linear transformation suffices for finding $\tilde{K}$.

The constants $D_\beta$ satisfy
%$D_\beta \simeq \bigl(C_{\frac{3\beta}{1-3\beta}}\bigr)^{1-3\beta} = O\Bigl(\bigl(\frac{\beta}{2-9\beta}\bigr)^{3\beta/2}\Bigr)$ 
$D_\beta \simeq \bigl(C_{\frac{4\beta}{2-3\beta}}\bigr)^{(2-3\beta)/2} = O\Bigl(\left(2-5\beta\right)^{-\beta}\Bigr)$
as $\beta \to \frac{2}{5}^-$, where $C_\alpha$ with $\alpha = \frac{4\beta}{2-3\beta}$ 
are the constants appearing in Theorem \ref{thm:Pisier-regular-ellipsoids}.
\end{theorem}

\medskip

As is clear from the statement, the quantitative estimates here are worse than in the symmetric case. Moreover, as will be clear from the proofs, we do not analyse and try to adapt Pisier's proof of Theorem \ref{thm:Pisier-regular-ellipsoids}, but rather we use it as a black box, and, just as in the case of Milman's theorem, we resort to other tools that allow us to compare covering numbers of a non-symmetric convex body with covering numbers of associated symmetric bodies. Thus we expect that the estimates are far from optimal, and that a more direct approach would be needed in order to obtain 
%comparable estimates
estimates considerably closer to those in Pisier's theorem.  

\medskip

We should also observe that Pisier obtains Theorem \ref{thm:Pisier-regular-ellipsoids} by proving the following stronger statement first: for every $\alpha\in (0,2)$ and every symmetric convex body $K$ in ${\mathbb R}^n$ we can find a linear image $\tilde{K}$ of $K$ such that, for every $1\ls l\ls n$,
\begin{equation}\label{eq:Pisier-Gelfand-numbers}
\max\{c_l(\tilde{K},B_2^n),\,c_l(\tilde{K}^\circ, B_2^n)\}\ls C_0C_\alpha^{1/\alpha}\left(\frac{n}{l}\right)^{1/\alpha},
\end{equation}
where $c_l(K,L)$ is called the \emph{Gelfand number} of the convex body $K$ with respect to the convex body $L$ and is defined as follows:
\begin{equation*}
c_l(K,L):= \inf\{r: \ \exists\,F\in G_{n,n-l+1} \ \hbox{such that}\ K\cap F\subseteq r(L\cap F)\}
\end{equation*}
(in \eqref{eq:Pisier-Gelfand-numbers}, $C_\alpha$ is the same constant as in the statement of Theorem \ref{thm:Pisier-regular-ellipsoids}). Once this is established, one can invoke a theorem by Carl \cite{Carl-1981} (see Section \ref{sec:prelim} for its statement), which connects the Gelfand numbers with the covering numbers, to derive Theorem \ref{thm:Pisier-regular-ellipsoids} as a corollary. We can recover the analogue of this stronger statement too.

\begin{proposition}\label{prop:Gelfand-numbers}
For every $\beta \in \bigl(0,\frac{2}{5}\bigr)$, and every convex body $K$ in ${\mathbb R}^n$ {\bf which has either barycentre or Santaló point at the origin}, we can find a linear image $\tilde{K}$ of $K$ such that, for all $1\ls l\ls n$,
\begin{equation*}
\max\{c_l(\tilde{K},B_2^n),\,c_l(\tilde{K}^\circ, B_2^n)\}\ls C_0D_\beta^{1/\beta}\left(\frac{n}{l}\right)^{1/\beta},
\end{equation*}
where $D_\beta$ is the same constant as in Theorem \ref{main-result} and $C_0$ is an absolute constant.
\end{proposition}

One of the tools needed for reducing the non-symmetric case to the symmetric one is a Blaschke-Santaló-type inequality for projections of convex bodies which have Santaló point at the origin. Recall that given a convex body $L$, the function $z\in {\rm int}(L)\mapsto |L||(L-z)^\circ|$ is strictly convex and has a unique point of minimum, called the Santaló point $s_L$ of $L$ (here $|\cdot|$, or sometimes $|\cdot|_n$, will denote $n$-dimensional Lebesgue volume). The celebrated Blashke-Santaló inequality (first proved by Blaschke \cite{Blaschke-1917} in ${\mathbb R}^2$ and ${\mathbb R}^3$, and then by Santaló \cite{Santalo-1949} in all dimensions; see also \cite{Petty-1985}, \cite{SaintRaymond-1981}, \cite{Meyer-Pajor-1990} and \cite{Bianchi-Kelly-2015} for further proofs of the inequality and/or its equality cases) states that, over convex bodies $L$ for which $s_L=0$, the volume product $|L||L^\circ|$ is uniquely maximised at ellipsoids.

Moreover, by the facts that (i) $(L^\circ)^\circ = L$ and (ii) if $L$ has barycentre at the origin, then $s_{L^\circ} = 0$ (another characterisation for the Santaló point of $L^\circ$), we can also conclude that $|L||L^\circ|\ls |B_2^n|^2$ for all bodies $L$ with barycentre at the origin.

We give similar approximate results for projections of such convex bodies.

\begin{theorem}\label{thm:BS-ineq-proj-sec}
Let $K$ be a convex body in ${\mathbb R}^n$ which has either barycentre or Santaló point at the origin. For every $1\ls l < n$, and for every subspace $F\in G_{n,l}$, we have
\begin{equation}
\bigl(|{\rm Proj}_F(K)|_l\cdot |K^\circ\cap F|_l\bigr)^{1/l}\ls C_0\frac{n}{l} |B_2^l|^{2/l}
\end{equation}
for some absolute constant $C_0$. 

Furthermore, up to the value of the constant $C_0$, the conclusion is optimal: one example is given by projections of the $n$-dimensional simplex.
\end{theorem}

\section{Preliminaries}\label{sec:prelim}

In the inequalities that we consider or establish in this paper, we will use the symbols $C, \tilde{C}, c,\tilde{c}$, and so on, to denote absolute constants which do not depend on any of the parameters (and certainly not on the dimension of the ambient space) unless specified otherwise by using a subscript. The value of these unspecified constants may be different in different occurrences of the symbol. We will sometimes also write $a\lesssim b$ or $b\gtrsim a$ if there is an absolute constant such that $a\ls Cb$. Finally, $a\simeq b$ means that both $a\lesssim b$ and $b\lesssim a$ hold.

\smallskip

A convex body $K$ in ${\mathbb R}^n$ is a convex, compact subset of ${\mathbb R}^n$ with non-empty interior. We write $\langle\cdot,\cdot\rangle$ for the standard dot product in ${\mathbb R}^n$, while $|A|$ (or sometimes $|A|_n$) stands for the $n$-dimensional Lebesgue volume. If $l\in \{1,2,\ldots,n-1\}$, by $G_{n,l}$ we denote the Grassmannian space of all $l$-dimensional (linear) subspaces of ${\mathbb R}^n$. If $F\in G_{n,l}$, then ${\rm Proj}_F(K)$ is the orthogonal projection of $K$ onto $F$, and $K\cap F$ is its section by $F$ (assuming that $K\cap F\neq \emptyset$). Unless specified otherwise, when we write 
\begin{equation*}
|{\rm Proj}_F(K)| \quad \hbox{or}\ \ |K\cap F|,
\end{equation*}
we mean the $l$-dimensional Lebesgue volume of these sets (volume within the subspace $F$). Recall that, if $0\in {\rm int}(K)$, then the polar set $(K\cap F)^\circ = (K\cap F)^{\circ,F}$ of the section $K\cap F$ within the subspace $F$ is well-defined and equal to the projection ${\rm Proj}_F(K^\circ)$.

We write $B_2^n$ for the unit Euclidean ball in ${\mathbb R}^n$, and $S^{n-1}$ for the unit Euclidean sphere, which is endowed with a (unique) rotationally invariant probablity measure $\sigma$ (the Haar probability measure). Given a convex body $K$ with $0\in {\rm int}(K)$, we define the \emph{inradius} $r(K)$ and \emph{circumradius} $R(K)$ of $K$ by
\begin{equation*}
r(K) :=\sup\{r>0: rB_2^n\subseteq K\} \quad\hbox{and}\ \ R(K):=\inf\{R>0: K\subseteq RB_2^n\}.
\end{equation*}
Since $K\subseteq L\,\Rightarrow\,L^\circ\subseteq K^\circ$, we have that $r(K) = \frac{1}{R(K^\circ)}$ and $R(K)= \frac{1}{r(K^\circ)}$. The barycentre $b(K)$ of $K$ is defined as follows:
\begin{equation*}
b(K) : = \frac{1}{|K|}\int_K x\,dx.
\end{equation*}
We define the support function $h_K$ of $K$ by
\begin{equation*}
x\in {\mathbb R}^n\mapsto h_K(x):=\max\{\langle y,x\rangle: y\in K\},
\end{equation*}
and if we also have $0\in {\rm int}(K)$, we define its Minkowski functional by
\begin{equation*}
x\in {\mathbb R}^n\mapsto \|x\|_K:=\inf\{t>0: x\in tK\}
\end{equation*}
(if $K$ is origin-symmetric, then $\|\cdot\|_K$ is a norm on ${\mathbb R}^n$). 
We can check that $\|x\|_K= h_{K^\circ}(x)$ and $h_K(x) = \|x\|_{K^\circ}$. The \emph{mean width} of $K$ is defined by
\begin{equation*}
M^\ast(K):= \int_{S^{n-1}} h_K(\theta)\,d\sigma(\theta), 
%\quad \hbox{and}\ \ M(K):=\int_{S^{n-1}}\|\theta\|_K\,d\sigma(\theta)
\end{equation*}
%to be the mean width and the mean norm of $K$ respectively (where $\sigma$ is the Haar probability measure on the unit Euclidean sphere $S^{n-1}$). 
and its \emph{mean norm} by
\begin{equation*}
M(K):=\int_{S^{n-1}}\|\theta\|_K\,d\sigma(\theta).
\end{equation*}
Given what we just said, we have $M^\ast(K)= M(K^\circ)$ and $M(K)= M^\ast(K^\circ)$.

\smallskip

For the covering numbers, the following volumetric bounds hold: 
%$\frac{|K|}{|L|} \ls N(K,L)$, and if $L$ is origin-symmetric as well, then
\begin{equation}\label{eq:covering-volumetric-bounds}
\frac{|K|}{|L|}\ls N(K,L) \ls \frac{|K+\tfrac{1}{2}(L\cap (-L))|}{|\tfrac{1}{2}(L\cap (-L))|} = \frac{|2K+(L\cap (-L))|}{|L\cap (-L)|}
%\ls \frac{|K+\tfrac{1}{2}L|}{|\tfrac{1}{2}L|} = \frac{|2K+L|}{|L|}.
\end{equation}
(we can even assume that $L$ has been translated first so that $L\cap (-L)$ has `large' volume, or even that $|L\cap (-L)| = \sup\limits_{x\in {\rm int}(L)}|(L-x)\cap (x-L)|$).

\subsection{``Local theory'': results about projections and sections}

Given a convex body $K$ in ${\mathbb R}^n$, its \emph{volume radius} is defined by ${\rm vrad}(K):= \left(\frac{|K|}{|B_2^n|}\right)^{1/n}$ (observe that it is $\simeq \sqrt{n}\,|K|^{1/n}$). If $F\in G_{n,l}$, then ${\rm vrad}({\rm Proj}_F(K))$ is understood to mean the volume radius of ${\rm Proj}_F(K)$ within the subspace $F$, thus $\left(\frac{|{\rm Proj}_F(K)|}{|B_2^l|}\right)^{1/l}$; similarly for ${\rm vrad}(K\cap F)$.

%Given a convex body $K$, its \emph{volume radius} ${\rm vrad}(K)$ is defined by ${\rm vrad}(K):= \left(\frac{|K|}{|B_2^n|}\right)^{1/n}$ where $n$ is the dimension of the affine hull of $K$
%%the smallest linear (or affine) subspace in which $K$ has non-empty interior 
%(observe that ${\rm vrad}(K) \simeq \sqrt{n}\,|K|^{1/n}$). Thus, if $F\in G_{n,l}$, then ${\rm vrad}({\rm Proj}_F(K)) = \left(\frac{|{\rm Proj}_F(K)|}{|B_2^l|}\right)^{1/l}$; similarly for ${\rm vrad}(K\cap F)$.

\smallskip

For any $1\ls l<n$, we set
\begin{equation*}
v_l(K) := \sup_{F\in G_{n,l}} {\rm vrad}({\rm Proj}_F(K)) \quad\hbox{and}\ \ v_l^-(K) := \inf_{F\in G_{n,l}} {\rm vrad}({\rm Proj}_F(K)),
\end{equation*}
and if $0\in {\rm int}(K)$, we also set
\begin{equation*}
w_l(K) := \sup_{F\in G_{n,l}} {\rm vrad}(K\cap F) \quad\hbox{and}\ \ w_l^-(K) := \inf_{F\in G_{n,l}} {\rm vrad}(K\cap F).
\end{equation*}

For any $k\gs 1$, we define the \emph{entropy number} $e_k(K,L)$ of $K$ by $L$ as follows:
\begin{equation*}
e_k(K,L) : = \inf\{R: N(K,RL)\ls 2^{k-1}\}.
\end{equation*}
Recall also the Gelfand number $c_k(K,L)$ that was already mentioned in the introduction: for any $k\in \{1,2,\ldots,n\}$:
\begin{equation*}
c_k(K,L):= \inf\{r: \ \exists\,F\in G_{n,n-k+1} \ \hbox{such that}\ K\cap F\subseteq r(L\cap F)\}
\end{equation*}
(for convenience in what follows, we can also set $c_k(K,L)=0$ if $k>n$). The following theorem is what we referred to as Carl's theorem in the introduction, and it shows that the entropy numbers are dominated in a certain sense by the Gelfand numbers (that is, in terms of any `reasonable' Lorentz norm).

\begin{theorem} {\rm (Carl, \cite{Carl-1981}, 1981)}\label{thm:Carl-entropy-Gelfand}
Suppose $K,L$ are origin-symmetric convex bodies in ${\mathbb R}^n$. For any $\alpha>0$ there exist constants $\rho_\alpha, \rho^\prime_\alpha > 0$ such that, for every integer $k\gs 1$, we will have
\begin{equation*}
\sup_{m=1,2,\ldots,k} m^\alpha e_m(K,L)\ls \rho_\alpha \sup_{m=1,2,\ldots,k} m^\alpha c_m(K,L)
\end{equation*} 
and
\begin{equation*}
\sum_{m=1}^k m^{-1+\alpha}e_m(K,L)\ls \rho^\prime_\alpha\sum_{m=1}^k m^{-1+\alpha}c_m(K,L)\,.
\end{equation*}
\end{theorem}

\subsection{Volume comparison results}

\noindent \emph{1. The Bourgain-Milman inequality.} In their seminal work \cite{Bourgain-Milman-1987}, Bourgain and Milman established an `asymptotic inverse' to the Blaschke-Santaló inequality: there exists an absolute constant $c_0>0$ (independent of $n$) such that, for all convex bodies $K$ in ${\mathbb R}^n$ with $0\in {\rm int}(K)$,
\begin{equation}\label{eq:Bourgain-Milman}
|K||K^\circ|\gs |K||(K-s(K))^\circ| \gs c_0^n|B_2^n|^2.
\end{equation}
%(note that we don't have to assume anything about the location of $s(K)$ or $b(K)$).
Combining this with the Blaschke-Santaló inequality, we can conclude that
\begin{equation*}
\big(|K||K^\circ|\bigr)^{1/n} \simeq |B_2^n|^{2/n} \simeq \frac{1}{n}
\end{equation*}
if $K$ is origin-symmetric, or if it has barycentre or Santaló point at the origin.

Moreover, if $K$ is origin-symmetric, we have that
\begin{equation*}
{\rm vrad}({\rm Proj}_F(K))\simeq \frac{1}{{\rm vrad}(K^\circ\cap F)}
\end{equation*}
for any $F\in G_{n,l}$, and hence also that $v_l(K)\simeq 1/w_l^-(K^\circ)$ and $v_l^-(K)\simeq 1/w_l(K^\circ)$. (We will see that, if $L$ is not origin-symmetric, these equivalences may no longer hold.)

\medskip

\noindent \emph{2. Volume of the difference body and of the symmetric convex hull.} In \cite{Rogers-Shephard-1957} Rogers and Shephard show that, for any convex body $K$ in ${\mathbb R}^n$, 
\begin{equation*}
|K-K| \ls \binom{2n}{n}|K|
\end{equation*}
with equality if and only if $K$ is an $n$-dimensional simplex (here $K-K:=\{x-y:x,y\in K\}$ is the \emph{difference body} of $K$). Moreover, in \cite{Rogers-Shephard-1958} they also prove that, if $0\in K$, then
\begin{equation}\label{eq:ineq-sym-conv}
|{\rm conv}(K,-K)|\ls 2^n|K|.
\end{equation}
It follows that ${\rm vrad}(K)\simeq {\rm vrad}(K-K)\simeq {\rm vrad}({\rm conv}(K,-K))$.

\medskip

\noindent \emph{3. The Milman-Pajor inequality: volume of the symmetric intersection.} In \cite{Milman-Pajor-2000} Milman and Pajor prove that, if $K$ is a centred convex body in ${\mathbb R}^n$ ($b(K)=0$), then
\begin{equation}\label{eq:Milman-Pajor-ineq}
|K\cap (-K)|\gs \frac{1}{2^n}|K|.
\end{equation} 
Note that the base of the exponential may not be optimal (see e.g. \cite{HSTV-2021} and \cite{CHMT-2022}), but also that for the arguments here an exact value is not needed, and instead all we require is an inequality of the form:
% there is an absolute constant $c>0$ such that
\begin{equation}\label{eq:vol-sym-intersec}
|K\cap (-K)| \gs c^n |K|
\end{equation}
for some absolute constant $c>0$.
%(in fact, Milman and Pajor show that we can take $c=\frac{1}{2}$). 

As e.g. Rudelson has remarked in \cite{Rudelson-2000b}, for this it suffices to combine the Blaschke-Santaló and Bourgain-Milman inequalities, as well as the relation $\bigl({\rm conv}(K,-K)\bigr)^\circ = K^\circ\cap (-K^\circ)$, which is valid for any convex body $K$ with $0\in {\rm int}(K)$. Rudelson applies this reasoning to obtain \eqref{eq:vol-sym-intersec} in the case when $s(K)=0$.

\medskip

\begin{remark}\label{rem:reduction-non-sym-M-position}
We can now quickly explain how Milman's theorem about the existence of an $M$-position can be extended to non-symmetric convex bodies. Assume that $K$ has barycentre or Santaló point at the origin, and that $\overline{K}:={\rm conv}(K,-K)$ is in $M$-position. In other words, assume that we know
%\begin{multline*}
%\max\big\{N({\rm conv}(K,-K), B_2^n),\,N(B_2^n,{\rm conv}(K,-K)),\\
%N(K^\circ\cap (-K^\circ),B_2^n),\,N(B_2^n,K^\circ\cap (-K^\circ))\bigr\}\ls \exp(Cn).
%\end{multline*}
\begin{equation*}
\max\big\{N(\overline{K}, B_2^n),\,N(B_2^n,\overline{K}),\,N(K^\circ\cap (-K^\circ),B_2^n),\,N(B_2^n,K^\circ\cap (-K^\circ))\bigr\}\ls \exp(Cn).
\end{equation*}
Then the covering numbers $N(K,B_2^n)$ and $N(B_2^n, K^\circ)$ are already controlled. For the covering number $N(B_2^n, K)$, we instead use that it is $\ls N(B_2^n, K\cap (-K))$. By the Milman-Pajor inequality (or Rudelson's variation), by the Rogers-Shephard inequality, and by the volumetric bounds \eqref{eq:covering-volumetric-bounds} for the covering numbers, we can check that
\begin{equation*}
N(B_2^n, K\cap (-K)) \ls \exp(C^\prime n) N(B_2^n,\overline{K})\ls \exp(C^{\prime\prime} n),
\end{equation*}
where we also crucially use the fact that $N(\overline{K}, B_2^n)\ls \exp(Cn)$ (and thus the Brunn-Minkowski inequality can be essentially reversed).
Similarly we bound the covering numbers $N(\overline{K^\circ}, B_2^n)$ and $N(K^\circ, B_2^n)$ by comparing them to the number $N(K^\circ\cap (-K^\circ),B_2^n)$.
\end{remark}

\medskip

\noindent \emph{4. Rudelson's result about sections of the difference body.} Note that, given $F\in G_{n,l}$, we have
\begin{multline*}
{\rm Proj}_F(K-K) = {\rm Proj}_F(K) - {\rm Proj}_F(K)\quad 
\\ \hbox{and}\ \ {\rm Proj}_F({\rm conv}(K,-K)) = {\rm conv}({\rm Proj}_F(K), -{\rm Proj}_F(K)),
\end{multline*}
thus the volumes of the projections of the difference body and of the symmetric convex hull of $K$ can be controlled well by the volume of the corresponding projection of $K$. In \cite{Rudelson-2000a} Rudelson suggests an analogue of these results for the sections of the difference body: he shows that
\begin{equation}\label{eq:Rudelson-section-result}
|(K-K)\cap F|^{1/l} \ls C\min\Bigl\{\sqrt{l\,},\frac{n}{l}\Bigr\} \max_{x\in {\mathbb R}^n}|K\cap (x+F)|^{1/l}.
\end{equation}
Note that we will only make use of the bound with the $\frac{n}{l}$ factor.

\medskip

\noindent \emph{5. Fradelizi's result about the maximal-volume section.} In \cite{Fradelizi-1997} Fradelizi proves that, for any convex body $K$ in ${\mathbb R}^n$, any $1\ls l < n$ and any $F\in G_{n,l}$,
\begin{equation}\label{eq:Fradelizi-section-result}
\max_{x\in {\mathbb R}^n}|K\cap (x+F)|\ls \left(\frac{n+1}{l+1}\right)^l |K\cap (b(K)+F)|.
\end{equation}
He also settles the equality cases.

\medskip

\noindent \emph{6. Volume product of a projection and its orthogonal section.} Rogers and Shephard \cite{Rogers-Shephard-1958} (see also \cite{Chakerian-1967}) proved that, for any origin-symmetric convex body $K$ in ${\mathbb R}^n$, any $1\ls l<n$ and any $F\in G_{n,l}$,
\begin{equation}\label{eq:RS-Spingarn}
|K| \ls |{\rm Proj}_F(K)|\cdot |K\cap F^\perp| \ls \binom{n}{l} |K|
\end{equation}
(in fact the upper bound is shown to hold even when $0\in {\rm int}(K)$). Spingarn \cite{Spingarn-1993} extended the lower bound to centred convex bodies too.

\subsection{Isotropic convex bodies}\label{subsec:properties-isotropic}

A convex body $K$ in ${\mathbb R}^n$ is called \emph{isotropic} if $|K|=1$, the barycentre of $K$ is at the origin, and its covariance matrix is a multiple of the identity, or equivalently, for every $\theta\in S^{n-1}$ we have
\begin{equation*}
\int_K\langle x, \theta\rangle^2\,dx = L_K^2.
\end{equation*}
%(where $\langle\cdot,\cdot\rangle$ is the standard dot product). 
The constant $L_K$ is called the isotropic constant of $K$. Note that every convex body $\tilde{K}$ in ${\mathbb R}^n$ has an affine image which is in isotropic position; this affine image is unique up to orthogonal transformations.

%The study of properties of the isotropic position has been at the core of a wealth of research within Asymptotic Convex Geometry ever since it was asked by Bourgain \cite{Bourgain-1987} whether all isotropic constants can be uniformly bounded from above.

It has been a central question in Asymptotic Convex Geometry, first asked by Bourgain \cite{Bourgain-1986}, whether all isotropic constants can be uniformly bounded from above (at the time Bourgain needed uniform lower bounds instead, which he showed hold true as $L_K\gs L_{B_2^n}= \frac{1}{\sqrt{n+2}\,|B_2^n|^{1/n}} \simeq 1 $ for all $K\subset {\mathbb R}^n$). 
The first non-trivial bounds were by Bourgain \cite{Bourgain-1991} ($L_K\lesssim \sqrt[4]{n}\log(n)$ for all $K\subset {\mathbb R}^n$) and by Klartag \cite{Klartag-2006} ($L_K\lesssim \sqrt[4]{n}$). Very crucial reductions of the question to other central problems in Asymptotic Convex Geometry, and recent huge breakthroughs in those directions by Yuansi Chen \cite{Chen-2021}, Klartag and Lehec \cite{Klartag-Lehec-2022}, Jambulapati, Lee and Vempala \cite{JLV-2022}, and even more recently by Klartag \cite{Klartag-2023}, have led to logarithmic bounds in the dimension $n$ for $L_K$ (the second to last paper, from 2022, gave $L_K\lesssim \log^{2.2226}(n)$ for every $K\subset {\mathbb R}^n$, while the most recent one gives $L_K\lesssim \sqrt{\log(n)}\,$).

\medskip

Standard properties of the isotropic position that we will need in the sequel are the following:
\begin{itemize}
\item if $K$ is an origin-symmetric isotropic convex body in ${\mathbb R}^n$, then $K\supset L_KB_2^n$. If $K$ is not symmetric but only centred, then this inclusion may no longer hold, but there still is an absolute constant $c_0<1$ such that $K\supset c_0L_KB_2^n$ (see e.g. \cite[Subsection 3.2.1]{BGVV-book}). In other words, $r(K) = \frac{1}{R(K^\circ)} \gs c_0L_K$ in the isotropic position.
\item For every $1\ls l\ls n-1$, and for every $F\in G_{n,l}$ we have (see e.g. \cite[Proposition 5.1.15]{BGVV-book}) that
\begin{equation*}
|K\cap F^\perp|^{1/l}\simeq \frac{L_{\pi_F({\bm 1}_K)}}{L_K},
\end{equation*}
where $\pi_F({\bm 1}_K)$ stands for the marginal distribution of the uniform distribution on $K$ with respect to the subspace $F$, and $L_{\pi_F({\bm 1}_K)}$ is its isotropic constant (refer to \cite{BGVV-book} for the history and details of how one extends the definition of the isotropic constant to log-concave distributions, and how the general upper and lower bounds can be shown to be the same).

As a consequence, by \eqref{eq:RS-Spingarn} we have that, for $F\in G_{n,l}$,
\begin{equation}\label{eq:bounds4proj}
|{\rm Proj}_F(K)|^{1/l} \gs \frac{1}{|K\cap F^\perp|^{1/l}}\gtrsim \frac{L_K}{L_{\pi_F({\bm 1}_K)}} \gtrsim \frac{L_K}{\sup_{M\subset {\mathbb R}^l}L_M}\,.
\end{equation}
\end{itemize}

%We recall that $r(K)\gs cL_K$. Moreover, for every $F\in G_{n,m}$, we have
%\begin{equation}\label{eq:bounds4proj}
%{\rm vrad}({\rm Proj}_F(K)) \gs \frac{1}{|B_2^m|^{1/m}}\frac{1}{|K\cap F^\perp|^{1/m}} \gtrsim \sqrt{m}\,\frac{L_K}{\sup_{M\subset {\mathbb R}^m} L_M}
%\end{equation}
%by the results by Rogers-Shephard and by Spingarn, and by standard properties of the isotropic position (that we recalled in Section \ref{sec:prelim}).

\subsection{A suggested approach towards Theorem \ref{main-result} by Klartag and Milman}\label{subsec:KM-approach}

Just as in the case of Milman's theorem about the $M$-position, and the reasoning in Remark \ref{rem:reduction-non-sym-M-position}, Klartag and Milman \cite{Klartag-VMilman-preprint, Klartag-VMilman-journal} suggest sandwiching a convex body $K$ that contains the origin in its interior by the symmetric convex bodies $\overline{K} := {\rm conv}(K,-K)$ and $\underline{K}:= K\cap (-K)$ (for which we can apply Pisier's theorem). As we will see, the key question their approach leads and reduces to is how well we can compare 
\begin{equation*} |{\rm Proj}_F(\overline{K})|^{1/l} \qquad \hbox{and} \qquad |{\rm Proj}_F(\underline{K})|^{1/l} \end{equation*}
for arbitrary subspaces $F$ in all admissible dimensions $l$.

Indeed, one could assume (without loss of generality)
that $\underline{K}$ is in $\alpha$-regular $M$-position for some $\alpha\in (0,2)$.
%(note that Pisier's theorem applies to $\underline{K}$ and $\overline{K}$).
If it happens that $\overline{K}$ is in $\alpha$-regular $M$-position as well, then we are done
because all the covering numbers
\begin{equation*} N(K, t B_2^n),N(K^\circ, t B_2^n), N(B_2^n, tK), N(B_2^n, tK^\circ)\end{equation*}
can be upper-bounded by corresponding covering numbers involving $\underline{K}$ or $\overline{K}$.
However what is more likely is that the $\alpha$-regular $M$-ellipsoid ${\cal E} = {\cal E}(\overline{K},\alpha)$
of $\overline{K}$ is not going to be $B_2^n$ but some other ellipsoid. In this
case, we should first control the covering numbers 
%$N({\cal E}, tB_2^n)$, $t\gs 1$
\begin{equation} \label{eq:incomplete-approach1}
N({\cal E}, tB_2^n) \quad \hbox{and} \ \  N(B_2^n, t{\cal E}) = N({\cal E}^\circ, tB_2^n),
\qquad t\gs 1. \end{equation}
(In fact, estimating the covering numbers $N({\cal E}^\circ, tB_2^n)$ is not needed for the application
to the covering numbers of $K$, 
%but we will describe below how to bound the volume of projections of 
%${\cal E}^\circ$ as well in order to highlight the subtle difficulties of the other case).
but we sketch below this case as well in order to clarify the subtleties that make the other case more difficult.)
By the following lemma this reduces to controlling the volume of arbitrary orthogonal projections of ${\cal E}$
and ${\cal E}^\circ$.

\begin{lemma} {\rm (see e.g. \cite[Remark 5.15]{Pisier-book})} \label{lem:ellipsoid-covering}
Let ${\cal Q} = T(B_2^n)$ be an ellipsoid in ${\mathbb R}^n$ with $T$ positive definite,
and let $\lambda_1\gs \lambda_2\gs \cdots\gs \lambda_n$ be the eigenvalues of $T$.
For every $k\gs 1$ set
\begin{equation*}\phi_k(T) := \sup_{1\ls l\ls n}\Bigl\{2^{-k/l} \cdot \Bigl[\prod_{1\ls j\ls l} \lambda_j\Bigr]^{1/l}\Bigr\}. \end{equation*}
Then for every $k\gs 1$ we have 
$\ \phi_k(T) \ls e_{k+1}({\cal Q}, B_2^n) \ls 6\phi_k(T)$.
%\begin{equation*} \phi_k(T) \ls e_{k+1}(Q, B_2^n) \ls 6\phi_k(T). \end{equation*}
\end{lemma}
\begin{remark} \label{rem:ellipsoid-covering}
(i) Since $\lambda_1,\lambda_2,\ldots,\lambda_n$ is a decreasing sequence, for any $l<n$ we will have that
\begin{equation*}
\Bigl[\prod_{1\ls j\ls l+1} \lambda_j\Bigr]^{1/(l+1)} \ls \Bigl[\prod_{1\ls j\ls l} \lambda_j\Bigr]^{1/l}.
\end{equation*}
Also $2^{-k/l}\simeq 1$ for any $l\gs k$. We thus see that $\phi_K(T) \simeq \sup_{l\ls k}\Bigl\{2^{-k/l} \cdot \Bigl[\prod_{1\ls j\ls l} \lambda_j\Bigr]^{1/l}\Bigr\}$.
%\begin{equation*}
%\phi_K(T) \simeq \sup_{l\ls k}\Bigl\{2^{-k/l} \cdot \Bigl[\prod_{1\ls j\ls l} \lambda_j\Bigr]^{1/l}\Bigr\}.
%\end{equation*}

\smallskip

(ii) If we know that there exist some positive constants $\gamma$ and $C(\gamma) \gs 1$ such that, for all $1\ls l\ls n$,
\begin{equation} \Bigl[\prod_{1\ls j\ls l} \lambda_j\Bigr]^{1/l} \ls C(\gamma)\cdot\Bigl(\frac{n}{l}\Bigr)^\gamma, \end{equation}
then the numbers $\phi_k(T)$, and hence also the entropy numbers $e_{k+1}({\cal Q}, B_2^n)$, can be estimated from above as follows:
\begin{equation*} \phi_k(T) \simeq e_{k+1}({\cal Q}, B_2^n) \lesssim 
\left\{\begin{array}{cl} C(\gamma)\left(\frac{\gamma}{e\log(2)}\right)^\gamma\left(\frac{n}{k}\right)^\gamma & \hbox{if}\ k\in \bigl[\frac{\gamma}{\log(2)},\, \frac{\gamma}{\log(2)}\, n\bigr] \smallskip \\[0.6em]
\frac{C(\gamma)\, n^\gamma}{2^k} & \hbox{if}\ k< \frac{\gamma}{\log(2)}  \smallskip \\[0.5em]
\frac{C(\gamma)}{2^{k/n}}  & \hbox{if}\ k > \frac{\gamma}{\log(2)}\, n
 \end{array} \right..  \end{equation*}
Then by the definition of the entropy numbers we get
 \begin{equation*} N({\cal Q}, tB_2^n) \ls \exp\left(\gamma[6C(\gamma)]^{1/\gamma}\, \frac{n}{t^{1/\gamma}}\right) \qquad \hbox{for every}\ t\gs 6C(\gamma)/e^\gamma. \end{equation*}

(iii) It can be shown using the Cauchy-Binet formula (see e.g. \cite[Lemma 4.1]{Klartag-VMilman-journal}) 
that, for any $1\ls l< n$,
\begin{equation*} \prod_{1\ls j\ls l} \lambda_j  = \frac{1}{|B_2^l|}\cdot \max_{F\in G_{n,l}} |{\rm Proj}_F({\cal Q})| 
= \frac{1}{|B_2^l|}\cdot \max_{F\in G_{n,l}} |{\cal Q}\cap F|, \end{equation*} 
and similarly
\begin{equation*} 
\prod_{n-l+1\ls j\ls n} \lambda_j = 
 \frac{1}{|B_2^l|}\cdot \min_{F\in G_{n,l}} |{\cal Q}\cap F| = |B_2^l|\cdot \frac{1}{\max\limits_{F\in G_{n,l}} |{\rm Proj}_F({\cal Q}^\circ)|}.
\end{equation*} 
\end{remark}

Lemma \ref{lem:ellipsoid-covering} and Remark \ref{rem:ellipsoid-covering}(iii) show that,
to upper-bound the covering numbers in \eqref{eq:incomplete-approach1},
we need to upper-bound the volume of projections of ${\cal E}$ and ${\cal E}^\circ$ in different dimensions.
%(in fact, estimating the covering numbers $N({\cal E}^\circ, tB_2^n)$ is not needed for the application
%to the covering numbers of $K$, 
%%but we will describe below how to bound the volume of projections of 
%%${\cal E}^\circ$ as well in order to highlight the subtle difficulties of the other case).
%but we sketch below this case as well in order to clarify the subtleties that make the other case more difficult).
Klartag and Milman use the assumptions made above:
if for simplicity we assume we have considered a $1$-regular $M$-position for $\underline{K}$
and a $1$-regular $M$-ellipsoid for $\overline{K}$, then we have
\begin{multline*}
\max\bigl\{N\bigl(\underline{K}, t B_2^n\bigr),N\bigl((\underline{K})^\circ, t B_2^n\bigr), N\bigl(B_2^n, t\underline{K}\bigr), 
N\bigl(B_2^n, t(\underline{K})^\circ\bigr)\bigr\}
\ls \exp\bigl(C_1\, n/t\bigr)
\\
\hbox{and}\ \  \max\bigl\{N\bigl(\overline{K}, t {\cal E}\bigr),N\bigl((\overline{K})^\circ, t {\cal E}^\circ\bigr), 
N\bigl({\cal E}, t\overline{K}\bigr), N\bigl({\cal E}^\circ, t(\overline{K})^\circ\bigr)\bigr\} \ls \exp\bigl(C_1\, n/t\bigr)
\end{multline*}
for every $t\gs C_1$.

It follows immediately that, for any $1\ls l < n$ and any $F\in G_{n,l}$,
\begin{equation*} N\bigl({\rm Proj}_F({\cal E}^\circ), t\,{\rm Proj}_F((\overline{K})^\circ)\bigr) \ls 
N\bigl({\cal E}^\circ, t(\overline{K})^\circ\bigr) \ls 
\exp\bigl(C_1\, n/t\bigr) \end{equation*}
for every $t\gs C_1$, and if we set $t_l= C_1n/l$ we see that
\begin{equation*} N\bigl({\rm Proj}_F({\cal E}^\circ), t_l\,{\rm Proj}_F((\overline{K})^\circ)\bigr) \ls e^l. \end{equation*}
Therefore
\begin{align*} 
|{\rm Proj}_F({\cal E}^\circ)|^{1/l} &\ls \left[N\bigl({\rm Proj}_F({\cal E}^\circ), t_l\,{\rm Proj}_F((\overline{K})^\circ)\bigr)\right]^{1/l}\cdot 
|t_l\,{\rm Proj}_F((\overline{K})^\circ)|^{1/l}
\\
& \ls e\,t_l\, |{\rm Proj}_F((\overline{K})^\circ)|^{1/l} = eC_1 \frac{n}{l}\, |{\rm Proj}_F((\overline{K})^\circ)|^{1/l}.
\end{align*}
Similarly we can show that
\begin{equation*} |{\rm Proj}_F((\underline{K})^\circ)|^{1/l} \ls eC_1 \frac{n}{l} \,|B_F|^{1/l} = eC_1 \frac{n}{l} \,|B_2^l|^{1/l}. \end{equation*}
We conclude that
\begin{align*} \max_{F\in G_{n,l}}|{\rm Proj}_F({\cal E}^\circ)|^{1/l}& \ls 
eC_1 \frac{n}{l}\, \max_{F\in G_{n,l}}|{\rm Proj}_F((\overline{K})^\circ)|^{1/l}
\\
&\ls eC_1 \frac{n}{l}\, \max_{F\in G_{n,l}}|{\rm Proj}_F((\underline{K})^\circ)|^{1/l} \ls (eC_1)^2 \left(\frac{n}{l}\right)^2\,|B_2^l|^{1/l} \end{align*}
given that $\underline{K} \subseteq \overline{K} \Rightarrow (\overline{K})^\circ\subseteq (\underline{K})^\circ$,
%since this holds for any $F\in G_{n,l}$,
%\begin{equation*} \max_{F\in G_{n,l}} |{\rm Proj}_F({\cal E}^\circ)|^{1/l} \lesssim \left(\frac{n}{l}\right)^2\,|B_2^l|^{1/l},\end{equation*}
which gives an upper bound of the form discussed in Remark \ref{rem:ellipsoid-covering}(ii).
Observe that the same upper bound holds even when $l=n$, that is,
we have $|{\cal E}^\circ|^{1/n} \ls (eC_1)^2\,|B_2^n|^{1/n}$ by the same reasoning.

\medskip

Obtaining a similar upper bound for the volume of projections of ${\cal E}$ is trickier 
since the inclusion $\underline{K} \subseteq \overline{K}$ does not work in our favour anymore.
%This is exactly the point the approach of Klartag and Milman was getting incorrectly. 
Let us describe how Klartag and Milman deal with this.
%attempt to compare
%\begin{equation*} |{\rm Proj}_F(\overline{K})|^{1/l} \qquad \hbox{and} \qquad |{\rm Proj}_F(\underline{K})|^{1/l} \end{equation*}
%in the non-trivial direction
%(which, in a very similar way to above, would then allow us to compare $|{\rm Proj}_F({\cal E})|^{1/l}$ to $|B_2^l|^{1/l}$). 

\begin{remark} 
Note that so far we have only prescribed that $K$ contain the origin in its interior (so that we can consider $K^\circ$).
We could choose a more specific position for $K$: what is assumed in \cite{Klartag-VMilman-preprint,Klartag-VMilman-journal},
which is \underline{no longer} without loss of generality, 
is that $K$ has both its barycentre and its Santaló point at the origin,
or equivalently that both $K$ and $K^\circ$ are centred. 

Of course not every convex body in ${\mathbb R}^n$ has this property.
%its barycentre and its Santaló point coinciding,
%and therefore not every convex body can have a position that satisfies the above.
In fact, as has been examined by Meyer, Schütt and Werner \cite{Meyer-Schutt-Werner-2011},
there are convex bodies in all dimensions whose barycentre and Santaló point are far apart.
\end{remark}

Recall that, if $K$ is centred, then, by the Rogers-Shephard and Milman-Pajor inequalities,
\begin{equation}\label{eq0:difficult-step}
 |\overline{K}|^{1/n} \ls 2|K|^{1/n} \ls 4|\underline{K}|^{1/n}. \end{equation}
For lower-dimensional projections Klartag and Milman argue as follows:
\begin{description}
\item[To upper-bound $\bm{|{\rm Proj}_F(\overline{K})|^{1/l}}$] By the Rogers-Shephard inequality again,
\begin{equation*} |{\rm Proj}_F(\overline{K})|^{1/l} = \big|\overline{{\rm Proj}_F(K)}\big|^{1/l} \ls 2 |{\rm Proj}_F(K)|^{1/l}.  \end{equation*}
Next we would like to compare $|{\rm Proj}_F(K)|^{1/l}$ to $|K^\circ\cap F|^{1/l}$, but this is not
immediate from the classical Blaschke-Santaló inequality given that
${\rm Proj}_F(K)$ may have neither barycentre nor Santaló point at the origin even though $K$ is centred.

Instead Klartag and Milman establish an approximate Blaschke-Santaló-type inequality for ${\rm Proj}_F(K)$ as follows:
${\rm Proj}_F(K)$ is the support of the marginal $\pi_F({\bm 1}_K)\equiv \pi_F(K)$ of ${\bm 1}_K$ on $F$,
which has a $1/(n-l)$-concave density given by
\begin{equation*} y\in F \quad \mapsto \quad 
g_{\pi_F(K)}(y):= |K\cap (y+F^\perp)| = \int_{F^\perp}{\bm 1}_K(y+z) \,dz.\end{equation*}
%By the Brunn-Minkowski inequality we know that $\pi_F(K)$ has a $1/(n-l)$-concave density,
It is also not difficult to see that $\pi_F(K)$ has barycentre at the origin. 
It is then known that we can find an associated convex body $T\in F$
that will be centred and not too different from the support of $\pi_F(K)$ in the following sense:
\begin{equation} \label{eq:inclusions-for-barycentre}
\frac{1}{e}T\subseteq {\rm Proj}_F(K) \subseteq e\frac{n+1}{l+1} T
\end{equation}
(this is done by picking the right candidate from the family of sets 
\begin{equation*} K_p(g) := \Bigl\{y\in F: \int_0^\infty g(sy)s^{p-1}\,ds\gs g(0)/p\Bigr\}, 
%K_p(g_{\pi_F(K)}) := \Bigl\{y\in F: \int_0^\infty g_{\pi_F(K)}(sy)s^{p-1}\,ds\gs g_{\pi_F(K)}(0)/p\Bigr\}, 
\qquad p>0,  \end{equation*}
linked to a log-concave density $g$ with $g(0)>0$; these were shown by Ball \cite{Ball-1988} to be convex bodies for all $p\gs 1$ (see also \cite{Barthe-thesis}, where it's verified that the proof extends, with only minor changes, to the not-necessarily-even-density case too, as well as to $p\in (0,1)$);
the right candidate here is $T = K_{l+1}(g_{\pi_F(K)})$, which is centred
if the $l$-dimensional density $g = g_{\pi_F(K)}$ is centred; 
the inclusions \eqref{eq:inclusions-for-barycentre} can be shown 
via a straightforward adaptation of \cite[Lemma 2.2]{Klartag-2005}
combined with Fradelizi's result \eqref{eq:Fradelizi-section-result}).

Once we have \eqref{eq:inclusions-for-barycentre} with some centred body $T\in F$,
we can write
\begin{align*}
|{\rm Proj}_F(K)|^{1/l} \cdot |K^\circ\cap F|^{1/l} &\ls \Big|e\frac{n+1}{l+1} T\Big|^{1/l}\cdot |eT^\circ|^{1/l}
\\
&= e^2\frac{n+1}{l+1}\, \bigl(|T||T^\circ|\bigr)^{1/l} \ls e^2\frac{n+1}{l+1}\,|B_2^l|^{2/l},
\end{align*}
which gives
\begin{equation}\label{eq1:difficult-step}
\frac{|{\rm Proj}_F(K)|^{1/l}}{|B_2^l|^{1/l}} \ls e^2\frac{n+1}{l+1}\, \frac{|B_2^l|^{1/l}}{|K^\circ\cap F|^{1/l}}.
\end{equation}

In the next section we are going to describe a somewhat different way to obtain \eqref{eq1:difficult-step},
which takes advantage of more recent results
(which however still rely on the $1/(n-l)$-concavity of the functions 
$y\in F \mapsto |K\cap (y+F^\perp)|$).

\smallskip

\item[To lower-bound $\bm{|{\rm Proj}_F(\underline{K})|^{1/l}}$] 
The goal here is again to relate $|{\rm Proj}_F(\underline{K})|^{1/l}$ to $|K^\circ\cap F|^{1/l}$,
%in a non-trivial direction, 
and Klartag and Milman do so by assuming that $K^\circ$ is centred too.

%Recall first that $(\underline{K})^\circ = \overline{(K^\circ)}$, and hence, by the Bourgain-Milman inequality,
%\begin{equation*}
%\frac{|{\rm Proj}_F(\underline{K})|^{1/l}}{|B_2^l|^{1/l}}\gs c_0\,\frac{|B_2^l|^{1/l}}{\big|\overline{(K^\circ)}\cap F\big|^{1/l}}
%\end{equation*}
%for some absolute constant $c_0>0$.
%
%Secondly, by Rudelson's result \cite{Rudelson-2000}, we have that
%\begin{equation*}
%\frac{1}{\big|\overline{(K^\circ)}\cap F\big|^{1/l}} \gs c_1\,\frac{l}{n} \,\frac{1}{\max\limits_{y\in F^\perp}|K^\circ\cap (y+F)|^{1/l}}
%\end{equation*}
%with $c_1>0$ another absolute constant.

By the Bourgain-Milman inequality and by Rudelson's result, we have that
\begin{equation*}
\frac{|{\rm Proj}_F(\underline{K})|^{1/l}}{|B_2^l|^{1/l}}\gs c_0\,\frac{|B_2^l|^{1/l}}{\big|\overline{(K^\circ)}\cap F\big|^{1/l}}
\gs c_0^\prime \,\frac{l}{n} \,\frac{|B_2^l|^{1/l}}{\max\limits_{y\in F^\perp}|K^\circ\cap (y+F)|^{1/l}}
\end{equation*}
for some absolute constants $c_0, c_0^\prime>0$.

Now, since $K^\circ$ has been assumed centred, by Fradelizi's result we also have
\begin{equation*}
\frac{1}{\max\limits_{y\in F^\perp}|K^\circ\cap (y+F)|^{1/l}} \gs \frac{l+1}{n+1}\,\frac{1}{|K^\circ\cap F|^{1/l}}.
\end{equation*}
We therefore conclude that 
\begin{equation} \label{eq2:difficult-step}
\frac{|{\rm Proj}_F(\underline{K})|^{1/l}}{|B_2^l|^{1/l}} \gs c\,\Bigl(\frac{l}{n}\Bigr)^2\frac{|B_2^l|^{1/l}}{|K^\circ\cap F\big|^{1/l}}
\end{equation}
for some absolute constant $c>0$,
which combined with \eqref{eq1:difficult-step} 
%and \eqref{eq2:difficult-step} 
finally gives
\begin{equation*} 
|{\rm Proj}_F(\overline{K})|^{1/l}\ls 2 |{\rm Proj}_F(K)|^{1/l} \ls C \,\Bigl(\frac{n}{l}\Bigr)^3\,|{\rm Proj}_F(\underline{K})|^{1/l}.
\end{equation*}
\end{description}

\smallskip

One can now check that the above reasoning establishes the existence 
of $\alpha$-regular ellipsoids in the sense of Theorem \ref{thm:Pisier-regular-ellipsoids} 
for all non-symmetric convex bodies $K$ with the property that their barycentre and Santaló point coincide
(with $\alpha$ up to $2/9$ non-inclusive, or even up to $2/5$ non-inclusive with an additional trick, as we will see via more careful computations below). 

For example, this establishes the existence of regular ellipsoids for the simplex,
although the regularity of the covering numbers in this instance can in fact be chosen much stronger; 
the latter may be seen by a more direct, and more specific to the simplex, argument
which we detail in the concluding remarks of this paper.
  
%Of course, as already mentioned, not every non-symmetric convex body has this property.
%%barycentre and Santaló point coinciding. 
In the next section we slightly modify the above approach to also handle bodies which do not have the property $b(K)=s(K)$. 
Clearly it is the last two steps of the argument we need to be careful about: 
as we will see, we will keep the second one as is, and 
only concern ourselves with the first one. In other words, we will assume that $K^\circ$ is centred,
or equivalently that $K$ has Santaló point at the origin, 
and will establish approximate Blaschke-Santaló-type inequalities for the projections of $K$ in this positioning too.

\section{Blaschke-Santaló-type inequalities for projections: proof of Theorem \ref{thm:BS-ineq-proj-sec}}

%The main result of this section is Proposition \ref{prop:BS-ineq-proj-sec}.
%To prove it, 
We apply crucially the key theorem of \cite{Meyer-Pajor-1990},
where Meyer and Pajor reprove the classical Blaschke-Santaló inequality
for all (not-necessarily symmetric) convex bodies via Steiner symmetrisation.  
We will use their terminology and say that, given a convex body $K$ in ${\mathbb R}^n$, 
an affine hyperplane $H = \{x\in {\mathbb R}^n: \langle x,u_H\rangle = a_H\}$ and $\lambda\in (0,1)$, 
$H$ is \emph{$\lambda$-separating} for $K$ if
%Given a convex body $K$ in ${\mathbb R}^n$, 
%an affine hyperplane $H = \{x\in {\mathbb R}^n: \langle x,u_H\rangle = a_H\}$ and $\lambda\in (0,1)$, 
%we will use their terminology and say that $H$ is \emph{$\lambda$-separating} for $K$ if
\begin{equation*} |K\cap H^+| := |\{x\in K: \langle x,u_H\rangle \gs a_H\}| = \lambda |K|. \end{equation*}
Note that, no matter what $\lambda\in (0,1)$ is, if $H$ is $\lambda$-separating for $K$, then ${\rm int}(K)\cap H\neq \emptyset$.

\begin{theorem} {\rm (Meyer-Pajor, \cite{Meyer-Pajor-1990})} \label{thm:Meyer-Pajor-BS}
Let $K$ be a convex body in ${\mathbb R}^n$, and $H$ an affine hyperplane 
that is $\lambda$-separating for $K$ for some $\lambda\in (0,1)$.
Then there exists $z\in {\rm int}(K)\cap H$ such that
\begin{equation}
|K|\, |(K-z)^\circ| \ls \frac{|B_2^n|^2}{4\lambda(1-\lambda)}.
\end{equation}

More precisely (see \cite[Lemma 2]{Meyer-Pajor-1990}),
the optimal choice of $z$ above (which minimises the volume product on the left-hand side) 
is the unique point $z_0 \in {\rm int}(K)\cap H$ 
with the property that the barycentre of $(K-z_0)^\circ$ lies in $H^\perp \equiv {\mathbb R} u_H$.
\end{theorem}

\smallskip

%\begin{proposition} \label{prop:BS-ineq-proj-sec}
%Let $K$ be a not-necessarily symmetric convex body in ${\mathbb R}^n$,
%and assume that either the Santaló point or the barycentre of $K$ is the origin.
%Then there exists 
%%$m\gs 1$ 
%%(in fact, applying a very recent result of... we can take $m=1$ in all cases)  
%%and 
%an absolute constant $C_0$ such that the following holds: 
%for every $1\ls l < n$ and every subspace $F\in G_{n,l}$, we have
%\begin{equation}\label{eq:BS-ineq-proj-sec}
%\bigl(|{\rm Proj}_F(K)|\cdot |K^\circ\cap F|\bigr)^{1/l} \ls C_0 \frac{n}{l}\,|B_2^l|^{2/l}.
%\end{equation}
%\end{proposition}
%\begin{remark}
%%If we fix the barycentre to be the origin, then we can take $m=1$.
%The conclusion when the barycentre is at the origin 
%follows already from the argument used by Klartag and Milman 
%which we sketched in the previous section. However we give one more argument below
%which is along very similar lines to the second part of the proof, 
%concerning the Santaló point, and which hopefully will help elucidate the key idea.
%
%%In the latter case (where, to the best of our knowledge, the conclusion is also new), 
%%we also have $m=1$ if we apply a very recent result by...
%%This result improves upon previous estimates 
%%%by Fradelizi, Meyer and Yaskin 
%%given in \cite{Fradelizi-Meyer-Yaskin-2017};
%%if we use these instead, we end up with $m=1$ again.  
%\end{remark}
\noindent {\it Proof of Theorem \ref{thm:BS-ineq-proj-sec}.}
The desired inequality when $b(K)=0$ has already been established by Klartag and Milman, in the way that we sketched in the previous section. However we can give one more argument
which is along very similar lines to the second part of the proof 
(concerning the Santaló point), and which hopefully will help elucidate the key idea. 

Stephen and Zhang \cite{Stephen-Zhang-2017} have shown that,
for every $1\ls l < n$, for every $F\in G_{n,l}$ and for every $\xi \in F\setminus\{0\}$,
if we set $\xi^+_F := \{y\in F: \langle y, \xi\rangle \gs 0\}$, we have
\begin{equation} \label{eq:Stephen-Zhang-2017}
|{\rm Proj}_F(K)\cap \xi^+_F|_l = |\{y\in {\rm Proj}_F(K): \langle y, \xi\rangle \gs 0\}|_l 
\gs \Bigl(\frac{l}{n+1}\Bigr)^l \,|{\rm Proj}_F(K)|_l
\end{equation}
(recall that we are comparing $l$-dimensional volumes, thus the subscript $l$, which we will suppress from here on).
Applying this result with $-\xi$ as well, we also see that
\begin{equation*} |{\rm Proj}_F(K)\cap \xi^+_F| \ls \left(1- \Bigl(\frac{l}{n+1}\Bigr)^l\right) \,|{\rm Proj}_F(K)|. \end{equation*} 
Let us rewrite this result as follows: for every $\xi \in F\setminus\{0\}$
the subspace $\xi^\perp_F:= F\cap \xi^\perp$, which is an affine hyperplane of $F$ passing through the origin,
is $\lambda$-separating for ${\rm Proj}_F(K)$ for some $\lambda$ satisfying
\begin{equation*}  \lambda \in \left[\Bigl(\frac{l}{n+1}\Bigr)^l, \ 1- \Bigl(\frac{l}{n+1}\Bigr)^l\right].\end{equation*}
%(it is not difficult to see that for $n\gs 2$ and $1\ls l<n$, $[l/(n+1)]^l < \frac{1}{2} < 1-[l/(n+1)]^l$).

The idea is to apply Theorem \ref{thm:Meyer-Pajor-BS} with the right hyperplane chosen from the above
(that is, with the right $\xi_0\in F$) in order to conclude that
\begin{equation} \label{eqp1:prop:BS-ineq-proj-sec}
|{\rm Proj}_F(K)|\,|K^\circ\cap F|\ls \frac{|B_2^l|^2}{4\lambda_{\xi_0}(1-\lambda_{\xi_0})} 
\ls \frac{1}{2} \Bigl(\frac{n+1}{l}\Bigr)^l\,|B_2^l|^2.
\end{equation}
To see that there is a suitable choice of $\xi_0\in F$ leading to this conclusion, note the following:
\begin{itemize}
\item if the barycentre of $K^\circ\cap F$ is at the origin, then the classical Blaschke-Santaló inequality applies immediately
and we have $|{\rm Proj}_F(K)|\,|K^\circ\cap F|\ls |B_2^l|^2$.
\item If instead $b\bigl(K^\circ\cap F\bigr)\neq 0$, denote it by $\xi_0$; clearly $\xi_0\in F\setminus\{0\}$
in this case. Applying Theorem \ref{thm:Meyer-Pajor-BS}, we see that there exists 
$z_0\in {\rm int}({\rm Proj}_F(K)) \cap \xi_0^\perp$ such that
\begin{equation*} |{\rm Proj}_F(K)|\cdot\big|({\rm Proj}_F(K)-z_0)^\circ\big| \ls \frac{|B_2^l|^2}{4\lambda_{\xi_0}(1-\lambda_{\xi_0})},\end{equation*}
where the polar set here is taken within $F$. In fact, the optimal $z_0$ is the \underline{unique} point 
$\in {\rm int}({\rm Proj}_F(K)) \cap \xi_0^\perp$ such that $({\rm Proj}_F(K)-z_0)^\circ$ has barycentre $\in {\mathbb R}\xi_0$.
By the fact that $({\rm Proj}_F(K)-0)^\circ = K^\circ\cap F $ has this property already, we conclude that 0 is this unique point
and that \eqref{eqp1:prop:BS-ineq-proj-sec} holds.
\end{itemize}

\bigskip

We now assume that $K$ has Santaló point at the origin, or in other words that $K^\circ$ has barycentre at the origin. 
Fix $1\ls l < n$ and $F\in G_{n,l}$.
Again, if we manage to say that, for every $\xi\in F\setminus \{0\}$, the hyperplane $\xi^\perp_F$
is $\lambda$-separating for $K^\circ\cap F$ for some $\lambda$ satisfying
\begin{equation*} \lambda \in [c_{n,l}, 1-c_{n,l}] \end{equation*}
with $c_{n,l}$ very similar to the previous lower bound, we can then argue completely analogously to above (with $K^\circ\cap F$ taking on the role of ${\rm Proj}_F(K)$ in the proof, and vice versa), and conclude that
\begin{equation} \label{eqp2:prop:BS-ineq-proj-sec}
|{\rm Proj}_F(K)|\,|K^\circ\cap F|\ls \frac{|B_2^l|^2}{4\lambda_{\xi_0}(1-\lambda_{\xi_0})} 
\ls \frac{1}{2c_{n,l}}\,|B_2^l|^2.
\end{equation}
We can do this because of an analogous result by Myroshnychenko, Stephen and Zhang \cite{MSZ-2018}, which states that, for $L\subset {\mathbb R}^n$ centred,
\begin{equation}\label{eq:MSZ-2018}
|L\cap \xi^+_F|_l = |\{y\in L\cap F: \langle y, \xi\rangle \gs 0\}|_l 
\gs \Bigl(\frac{l}{n+1}\Bigr)^l \,|L\cap F|_l
\end{equation}
(just as in \cite{Stephen-Zhang-2017}, this is also shown to be optimal). We apply this with $L=K^\circ$ (this also shows why we focused on $K^\circ\cap F$ here, and not on ${\rm Proj}_F(K)$).

\smallskip

{\small {\bf Side Remark.} It should be noted that our desired inequality is rather crude compared to these results, since in the end we aim to bound volume radii, and therefore additional constant factors raised to the power $l$ in an inequality such as \eqref{eq:MSZ-2018} would essentially not affect this application of it. Thus, even an earlier result towards the same goal as in \cite{MSZ-2018} would work here: in \cite[Corollary 9]{Fradelizi-Meyer-Yaskin-2017} Fradelizi, Meyer and Yaskin show that 
\begin{equation*} c_{n,l} \gs c\Bigl(\frac{n+1}{n-l+1}\Bigr)^2\,\frac{1}{l^2}\,\Bigl(\frac{l}{n+1}\Bigr)^l > \frac{c}{l^2}\,\Bigl(\frac{l}{n+1}\Bigr)^l \end{equation*}
for some absolute constant $c>0$, which is $\gs \bigl(c^\prime\frac{l}{n+1}\bigr)^l$.}

%recall estimates by Fradelizi, Meyer and Yaskin \cite{Fradelizi-Meyer-Yaskin-2017} 
%or by Myroshnychenko, Stephen and Zhang \cite{MSZ-2018} (the latter, very recent work replacing the estimates in \cite{Fradelizi-Meyer-Yaskin-2017}
%by optimal ones, and also settling equality cases). Note that for our application, where multiplying by an
%absolute constant raised to the dimension of the ambient space makes no real difference, 
%the two results lead actually to the same conclusion. 
%
%Let us verify this by considering the estimates from \cite{Fradelizi-Meyer-Yaskin-2017}:
%applied to $K^\circ$ which is centred, these give
%\begin{equation*} c_{n,l} > \frac{c}{l^2}\,\Bigl(\frac{l}{n+1}\Bigr)^l \gs \bigl(c^\prime\frac{l}{n+1}\bigr)^l \end{equation*}
%for some absolute constants $c,c^\prime >0$.
%Indeed, \cite[Corollary 9]{Fradelizi-Meyer-Yaskin-2017} implies that
%\begin{align*}
%|K^\circ\cap \xi^+_F| &\gs \frac{c}{(n-l+1)^2}\,\Bigl(\frac{l}{n+1}\Bigr)^{l-2}\,|K^\circ\cap F|
%\\
%\intertext{for every subspace $F\in G_{n,l}$ and every $\xi \in F\setminus\{0\}$,
%where the right-hand side can also be rearranged as}
%&\gs c\Bigl(\frac{n+1}{n-l+1}\Bigr)^2\,\frac{1}{l^2}\,\Bigl(\frac{l}{n+1}\Bigr)^l \,|K^\circ\cap F| > 
%\frac{c}{l^2}\,\Bigl(\frac{l}{n+1}\Bigr)^l  |K^\circ\cap F|.
%\end{align*}
%
%In contrast, \cite{MSZ-2018} leads to the more precise $c_{n,l} \gs \bigl(l/(n+1)\bigr)^l$
%(with equality being possible now).
%
%With the above estimates for $c_{n,l}$ plugged in \eqref{eqp2:prop:BS-ineq-proj-sec},
%the proof of the inequality is complete.

\bigskip
\medskip

We now discuss the optimality of Theorem \ref{thm:BS-ineq-proj-sec}. Let $S_n$ be a regular simplex of edge-length $\sqrt{2}$ which has barycentre at the origin. Note that, if we consider $l$ of its vertices as well as the average of the remaining $n+1-l$, the affine hull of these points coincides with an $l$-dimensional (linear) subspace $F_l$ (which thus contains the barycentre of $S_n$). Moreover, by either embedding $S_n$ in ${\mathbb R}^{n+1}$ in the standard way:
\begin{equation*}
S_n = {\rm conv}\Bigl\{e_i-\frac{1}{n+1}\bm{1}:1\ls i\ls n+1\Bigr\}
\end{equation*}
(where $\bm{1}=e_1+e_2+\ldots+e_{n+1}$), or by explicitly writing down its vertices in ${\mathbb R}^n$ (starting with $n$ of its vertices being the vectors $e_i$, and then translating properly so that the barycentre becomes the origin), we can directly calculate that
\begin{equation*}
|S_n\cap F_l| = \frac{\sqrt{n+1}}{l!\,\sqrt{n+1-l}}.
\end{equation*}

At the same time, $S_n^\circ= -(n+1)S_n$, and thus
\begin{equation*}
{\rm Proj}_{F_l}(S_n^\circ) = -(n+1){\rm Proj}_{F_l}(S_n).
\end{equation*}
We can thus write
\begin{align*}
\bigl(|{\rm Proj}_{F_l}(S_n^\circ)|\cdot |S_n\cap F_l|\bigr)^{1/l} &= (n+1)\bigl(|{\rm Proj}_{F_l}(S_n)|\cdot |S_n\cap F_l|\bigr)^{1/l}
\\
&\gs (n+1)|S_n\cap F_l|^{2/l} \simeq \frac{n+1}{l^2} \simeq \frac{n+1}{l} |B_2^l|^{2/l},
\end{align*}
which shows optimality of the Blaschke-Santaló-type inequalities we got, both in the case of a convex body with Santaló point at the origin, and in the case of a centred convex body.
%\end{proof}
\qed

\bigskip

In \cite{Dirksen-2017} Dirksen states the conjecture that, out of all sections of the regular simplex that pass through its barycentre, the volume of the section that we considered is maximum. He also establishes an upper bound for the volume of such sections which is asymptotically equivalent: if $E\in G_{n,l}$, then
\begin{equation*}
|S_n\cap E| \ls \frac{(\sqrt{l+1})^{\frac{l+1}{n+1}}}{l!}.
\end{equation*}
A less precise but similar bound follows from the above result and reasoning: indeed,
\begin{equation*}
\frac{1}{2}\left(\frac{n+1}{l}\right)^l |B_2^l|^2 \gs |{\rm Proj}_E(S_n^\circ)|\cdot |S_n\cap E| \gs (n+1)^l\, |S_n\cap E|^2.
\end{equation*}

In the case of the subspace $F_l$ that we considered above, we could even estimate the volume product $|{\rm Proj}_{F_l}(S_n^\circ)|\cdot |S_n\cap F_l|$ directly and more accurately (and only later compare it with the given upper bound). Indeed, by the Rogers-Shephard inequality
\begin{equation*}
|{\rm Proj}_{F_l}(S_n)|^{1/l} \gs \frac{1}{2}|{\rm Proj}_{F_l}(\overline{S_n})|
\end{equation*}
and we can recall that $\overline{S_n}$ can be viewed as a projection of the $\ell_1$-unit ball $B_1^{n+1}$:
\begin{equation*}
\overline{S_n} = {\rm conv}\Bigl\{\pm\bigl(e_i-\tfrac{1}{n+1}\bm{1}\bigr):1\ls i\ls n+1\Bigr\} = {\rm Proj}_{H_{\bm 1}}\bigl(B_1^{n+1}\bigr)
\end{equation*}
with $H_{\bm 1}$ the hyperplane orthogonal to the vector ${\bm 1}$. But if we suppose that
$F_l$ is the subspace spanned by the first $l$ vertices of $S_n$, the vectors $e_j-\frac{1}{n+1}{\bm 1}$, $1\ls j\ls l$, then it can be checked that 
\begin{equation*}
{\rm Proj}_{F_l}(\overline{S_n}) =  {\rm Proj}_{F_l}\bigl(B_1^{n+1}\bigr) = {\rm conv}\Bigl\{\pm\bigl(e_j-\tfrac{1}{n+1}\bm{1}\bigr):1\ls j\ls l\Bigr\}.
\end{equation*}
At the same time, we have that the average of the other vertices of $S_n$ is in the convex hull of the vectors $-\bigl(e_j-\frac{1}{n+1}{\bm 1}\bigr)$, $1\ls j\ls l$, and the origin, and thus
\begin{equation*}
{\rm conv}(S_n\cap F_l, -(S_n\cap F_l)) = {\rm conv}\Bigl\{\pm\bigl(e_j-\tfrac{1}{n+1}\bm{1}\bigr):1\ls j\ls l\Bigr\} = {\rm Proj}_{F_l}(\overline{S_n}).
\end{equation*}
This implies that
\begin{equation*}
(n+1)^l|S_n\cap F_l|^2 \ls |{\rm Proj}_{F_l}(S_n^\circ)|\cdot |S_n\cap F_l| \ls (2(n+1))^l|S_n\cap F_l|^2.
\end{equation*}

This also shows that the `extreme' case that we used to show optimality of Theorem \ref{thm:BS-ineq-proj-sec} is very far from any extreme cases in Fradelizi's result (an instance of the latter is when the considered $l$-dimensional subspace is parallel to the affine hull of an $l$-dimensional face of the simplex). It is unclear to us whether there is a way to better combine the two results.

\section{Proving Theorem \ref{main-result} and Proposition \ref{prop:Gelfand-numbers}}

%\noindent{\it Proof of Theorem \ref{main-result}.}
We begin with the 
%proof of the theorem.
\medskip\\
\noindent{\it Proof of the theorem.}
Most of the steps we need have been discussed or hinted at in Subsection \ref{subsec:KM-approach}.
We quickly go over the necessary adjustments to prove the result in the general case
and for all $\beta\in \bigl(0,\frac{2}{5}\bigr)$ (note also that a direct adjustment of the proof scheme already discussed would lead to a range of $\bigl(0,\frac{2}{9}\bigr)$ for $\beta$; instead we also insert a `balancing' step in the middle or the proof which will allow us to slightly improve the final estimates; we thank Alexander Litvak for suggesting this trick). 

\smallskip

%We will also use a `bootstrapping' step that helps improve the final estimates on regularity compared to what a direct adaptation of subsection \ref{subsec:KM-approach} would give; we thank Alexander Litvak for suggesting this trick. 

%\smallskip

Fix $\beta\in \bigl(0,\frac{2}{5}\bigr)$ and set $\alpha = \frac{4\beta}{2-3\beta}$;
note that $\alpha\in (0,2)$. Let also $K\subset {\mathbb R}^n$ be a convex body. 

\smallskip

{\bf Step 1.} This time we start by assuming that $K$ has Santaló point at the origin, 
and that $\underline{K} := K\cap (-K)$ is in $\alpha$-regular $M$-position, that is,
\begin{equation*}
\max\bigl\{N\bigl(\underline{K}, t B_2^n\bigr),N\bigl((\underline{K})^\circ, t B_2^n\bigr), 
N\bigl(B_2^n, t\underline{K}\bigr), N\bigl(B_2^n, t(\underline{K})^\circ\bigr)\bigr\}
\ls \exp\bigl(C_\alpha\, n/t^\alpha\bigr)
\end{equation*}
for every $t\gs C_\alpha^{1/\alpha}$. We also find an $\alpha$-regular $M$-ellipsoid ${\cal E}$
for $\overline{K} : = {\rm conv}(K,-K)$:
\begin{equation*}
\max\bigl\{N\bigl(\overline{K}, t{\cal E}\bigr), N\bigl((\overline{K})^\circ, t{\cal E}^\circ\bigr), 
N\bigl({\cal E}, t\overline{K}\bigr), N\bigl({\cal E}^\circ, t(\overline{K})^\circ\bigr)\bigr\} 
\ls \exp\bigl(C_\alpha\, n/t^\alpha\bigr)
\end{equation*}
for every $t\gs C_\alpha^{1/\alpha}$. 

Finally, along with the above assumptions we can make one more, which will help make Step 2 easier to read: by applying an orthogonal transformation to $K$ (and hence to $\overline{K}$ and $\underline{K}$ too) if needed, we can assume that ${\cal E} = \Delta_\lambda B_2^n$ with $\Delta_\lambda$ a diagonal matrix with diagonal entries $\lambda_1 \gs \lambda_2\gs \cdots\gs \lambda_n>0$.

%Then, for every $1\ls l < n$, we can obtain as before that
%\begin{equation*}
%\max_{F\in G_{n,l}}|{\rm Proj}_F({\cal E}^\circ)|^{1/l} \ls e^2C_\alpha^{2/\alpha} \left(\frac{n}{l}\right)^{2/\alpha}\,|B_2^l|^{1/l}.
%\end{equation*}
%Similarly $|{\cal E}^\circ|^{1/n} \ls e^2C_\alpha^{2/\alpha}\,|B_2^n|^{1/n}$.

\medskip

%For the remaining volume inequalities, 
We first observe that, by the Rogers-Shephard, Blaschke-Santaló and Bourgain-Milman inequalities,
\begin{align*} 
%|K^\circ|\ls |{\rm conv}(K^\circ, -K^\circ)| = |(\underline{K})^\circ|
\qquad &|(\underline{K})^\circ| = |{\rm conv}(K^\circ, -K^\circ)| \ls 2^n|K^\circ| \\ 
\Longrightarrow \quad  & |\overline{K}| \ls 2^n|K| \ls 4^n \frac{|K|\cdot |K^\circ|}{|(\underline{K})^\circ| \cdot |\underline{K}|}\ |\underline{K}| \ls C_0^n |\underline{K}| \\
\Longrightarrow \quad & |{\cal E}|^{1/n} \ls e^2C_0\, C_\alpha^{2/\alpha} |B_2^n|^{1/n}
\end{align*} 
for some absolute constant $C_0$ (depending on the constant in the Bourgain-Milman inequality).

Moreover, obtaining \eqref{eq1:difficult-step} by Theorem \ref{thm:BS-ineq-proj-sec}
(which we apply in the case that the Santaló point is at the origin), 
and combining this with \eqref{eq2:difficult-step} (which is true when $K^\circ$ is centred, as is the case here),
we can write, for every $1\ls l < n$ and every $F\in G_{n,l}$,
\begin{equation}\label{eqp0:thm:main-result} 
|{\rm Proj}_F(\overline{K})|^{1/l}\ls 2 |{\rm Proj}_F(K)|^{1/l} \ls C \,\Bigl(\frac{n}{l}\Bigr)^3\,|{\rm Proj}_F(\underline{K})|^{1/l}. 
\end{equation} 
We also have
\begin{equation*} 
|{\rm Proj}_F({\cal E})|^{1/l} \ls e\,C_\alpha^{1/\alpha}\Bigl(\frac{n}{l}\Bigr)^{1/\alpha}\, |{\rm Proj}_F(\overline{K})|^{1/l} 
\ \ \hbox{and}
\ \ |{\rm Proj}_F(\underline{K})|^{1/l} \ls e\,C_\alpha^{1/\alpha}\Bigl(\frac{n}{l}\Bigr)^{1/\alpha}\, |B_2^l|^{1/l}.
\end{equation*}
Therefore, 
%\begin{equation*}\label{eqp1:thm:main-result}
%\max_{F\in G_{n,l}}|{\rm Proj}_F({\cal E})|^{1/l} 
%\ls \tilde{C}\,C_\alpha^{2/\alpha} \left(\frac{n}{l}\right)^{(2+3\alpha)/\alpha}\,|B_2^l|^{1/l}.
%\end{equation*}
\begin{equation*}
v_l({\cal E}) = \max_{F\in G_{n,l}}\frac{|{\rm Proj}_F({\cal E})|^{1/l}}{|B_2^l|^{1/l}}\ls \tilde{C}\,C_\alpha^{2/\alpha} \left(\frac{n}{l}\right)^{(2+3\alpha)/\alpha}.
\end{equation*}

\medskip

{\bf Step 2.} Set $\Delta_{\!\sqrt{\lambda}}$ to be the diagonal matrix with diagonal entries $\sqrt{\lambda_1}\gs \sqrt{\lambda_2}\gs\cdots\sqrt{\lambda_n}$, and ${\cal Q}_{{\cal E}}$ to be the ellipsoid $\Delta_{\!\sqrt{\lambda}}B_2^n = \Delta_{\!\sqrt{\lambda}}^{-1}{\cal E}$. We have that
\begin{equation*}
\max\bigl\{N\bigl(\Delta_{\!\sqrt{\lambda}}^{-1}\underline{K}, t {\cal Q}_{{\cal E}}^\circ\bigr),N\bigl((\Delta_{\!\sqrt{\lambda}}^{-1}\underline{K})^\circ, t {\cal Q}_{{\cal E}}\bigr), 
N\bigl({\cal Q}_{{\cal E}}^\circ, t\Delta_{\!\sqrt{\lambda}}^{-1}\underline{K}\bigr), N\bigl({\cal Q}_{{\cal E}}, t(\Delta_{\!\sqrt{\lambda}}^{-1}\underline{K})^\circ\bigr)\bigr\}
\ls \exp\bigl(C_\alpha\, n/t^\alpha\bigr)
\end{equation*}
and 
\begin{equation*}
\max\bigl\{N\bigl(\Delta_{\!\sqrt{\lambda}}^{-1}\overline{K}, t{\cal Q}_{{\cal E}}\bigr), N\bigl((\Delta_{\!\sqrt{\lambda}}^{-1}\overline{K})^\circ, t{\cal Q}_{{\cal E}}^\circ\bigr), 
N\bigl({\cal Q}_{{\cal E}}, t\Delta_{\!\sqrt{\lambda}}^{-1}\overline{K}\bigr), N\bigl({\cal Q}_{{\cal E}}^\circ, t(\Delta_{\!\sqrt{\lambda}}^{-1}\overline{K})^\circ\bigr)\bigr\} 
\ls \exp\bigl(C_\alpha\, n/t^\alpha\bigr)
\end{equation*}
Furthermore, for all $1\ls l\ls n$,
%\begin{equation*}\label{eqp2:thm:main-result}
%\max_{F\in G_{n,l}}\frac{|{\rm Proj}_F({\cal Q}_{{\cal E}})|^{1/l}}{|B_2^l|^{1/l}} =  \sqrt{\max_{F\in G_{n,l}}\frac{|{\rm Proj}_F({\cal E})|^{1/l}}{|B_2^l|^{1/l}}\,}
%\ls \tilde{C}\,C_\alpha^{1/\alpha} \left(\frac{n}{l}\right)^{\frac{1}{\alpha} + \frac{3}{2}}.
%\end{equation*}
\begin{equation*}\label{eqp2:thm:main-result}
v_l({\cal Q}_{{\cal E}}) =  \sqrt{v_l({\cal E})}
\ls \tilde{C}\,C_\alpha^{1/\alpha} \left(\frac{n}{l}\right)^{\frac{1}{\alpha} + \frac{3}{2}}.
\end{equation*}

Now, applying Lemma \ref{lem:ellipsoid-covering} (as in Remark \ref{rem:ellipsoid-covering}(ii)), we see that
\begin{equation*} 
%\max\{N({\cal E}, tB_2^n), N({\cal E}^\circ, tB_2^n)\} 
N({\cal Q}_{{\cal E}}, tB_2^n) = N(B_2^n, t{\cal Q}_{{\cal E}}^\circ)
\ls \exp\left(\alpha^{-1}\tilde{C}_0\,C_\alpha^{2/(2+3\alpha)}\,\frac{n}{t^{\frac{2\alpha}{2+3\alpha}}}\right)
\end{equation*}
for every $t\gs \tilde{C}_0\,C_\alpha^{1/\alpha}$ for some absolute constant $\tilde{C}_0 \gs 2$.

{\bf Step 3.} Using the above, we can write
\begin{gather*}
%N(\overline{K}, tB_2^n)\ls 
N\bigl(\Delta_{\!\sqrt{\lambda}}^{-1}\overline{K}, tB_2^n\bigr) \ls N\bigl(\Delta_{\!\sqrt{\lambda}}^{-1}\overline{K}, s{\cal Q}_{{\cal E}}\bigr) N(s{\cal Q}_{{\cal E}}, tB_2^n)
 \ls \exp\bigl(C_\alpha\, n/s^\alpha\bigr) \cdot \exp\biggl(\alpha^{-1}\tilde{C}_0\,C_\alpha^{\frac{2}{2+3\alpha}}\,n \Bigl(\frac{s}{t}\Bigr)^{\frac{2\alpha}{2+3\alpha}}\biggr)
 \\ \nonumber
\hbox{as long as $\,s\gs C_\alpha^{1/\alpha}\ $ and $\ \dfrac{t}{s} \gs \tilde{C}_0\,C_\alpha^{1/\alpha}$.}
\end{gather*}
%\begin{equation*} s\gs C_\alpha^{1/\alpha}\qquad \hbox{and} \qquad  \frac{t}{s} \gs C_0\,C_\alpha^{2/\alpha}\end{equation*}
Optimising in $s$, we see that
\begin{equation*}\label{eq:opt-in-s}  \frac{C_\alpha}{s^\alpha} =  \frac{\alpha^{-1}\tilde{C}_0\, C_\alpha^{\frac{2}{2+3\alpha}}\,s^{\frac{2\alpha}{2+3\alpha}}}{t^{\frac{2\alpha}{2+3\alpha}}} \quad \Longrightarrow \quad s = \frac{C_\alpha^{\frac{3}{4+3\alpha}}\, t^{\frac{2}{4+3\alpha}}}{(\alpha^{-1}\tilde{C}_0)^{\frac{2+3\alpha}{4\alpha+3\alpha^2}}}, \end{equation*}
which is $\gs C_\alpha^{1/\alpha}$ if $t\gs (\alpha^{-1}\tilde{C}_0)^{\frac{2+3\alpha}{2\alpha}}\,C_\alpha^{2/\alpha}$.
Moreover, $\frac{t}{s} \gs \tilde{C}_0\,C_\alpha^{1/\alpha}$ 
%if $t\gs C_0^{\frac{3\alpha^2-2}{3\alpha^2+2\alpha}}\, C_\alpha^{3/\alpha}$
for the same $t$. Therefore
\begin{equation*}\label{eq:covering-num-1}
N\bigl(\Delta_{\!\sqrt{\lambda}}^{-1}\overline{K}, tB_2^n\bigr) \ls \exp\bigl(2C_\alpha\, n/s^\alpha\bigr) 
= \exp\Bigl(2(\alpha^{-1}\tilde{C}_0)^{\frac{2+3\alpha}{4+3\alpha}}\, C_\alpha^{\frac{4}{4+3\alpha}}\,\frac{n}{t^{\frac{2\alpha}{4+3\alpha}}}\Bigr) 
\end{equation*}
for every $t\gs (\alpha^{-1}\tilde{C}_0)^{\frac{2+3\alpha}{2\alpha}}\,C_\alpha^{2/\alpha}$.

\smallskip

Analogously, 
\begin{align*}
\nonumber
N\bigl(B_2^n, t\Delta_{\!\sqrt{\lambda}}^{-1}\underline{K}\bigr) &\ls N\Bigl(B_2^n, \frac{t}{s}{\cal Q}_{{\cal E}}^\circ\Bigr) 
N\Bigl(\frac{t}{s}{\cal Q}_{{\cal E}}^\circ, t\Delta_{\!\sqrt{\lambda}}^{-1}\underline{K}\Bigr)
\\ \nonumber 
&= N\Bigl({\cal Q}_{{\cal E}}, \frac{t}{s}B_2^n, \Bigr) 
N\Bigl({\cal Q}_{{\cal E}}^\circ, s\Delta_{\!\sqrt{\lambda}}^{-1}\underline{K}\Bigr)
\\ \nonumber
&\ls \exp\biggl(\alpha^{-1}\tilde{C}_0\,C_\alpha^{\frac{2}{2+3\alpha}}\,n \Bigl(\frac{s}{t}\Bigr)^{\frac{2\alpha}{2+3\alpha}}\biggr)
\cdot \exp\bigl(C_\alpha\, n/s^\alpha\bigr)
\\ \label{eq:covering-num-2}
& = \exp\Bigl(2(\alpha^{-1}\tilde{C}_0)^{\frac{2+3\alpha}{4+3\alpha}}\, C_\alpha^{\frac{4}{4+3\alpha}}\,\frac{n}{t^{\frac{2\alpha}{4+3\alpha}}}\Bigr) 
\end{align*}
with $s$ as before and $t \gs (\alpha^{-1}\tilde{C}_0)^{\frac{2+3\alpha}{2\alpha}}\,C_\alpha^{2/\alpha}$.

Finally, in the same way we can bound the covering numbers 
\begin{multline*}
N\bigl(B_2^n, t(\Delta_{\!\sqrt{\lambda}}^{-1}\overline{K})^\circ\bigr) \ls N\Bigl(B_2^n, \frac{t}{s}{\cal Q}_{{\cal E}}^\circ\Bigr) 
N\Bigl(\frac{t}{s}{\cal Q}_{{\cal E}}^\circ, t(\Delta_{\!\sqrt{\lambda}}^{-1}\overline{K})^\circ\Bigr)
\\
\hbox{and}\ \  N\bigl((\Delta_{\!\sqrt{\lambda}}^{-1}\underline{K})^\circ, tB_2^n\bigr) \ls N\bigl((\Delta_{\!\sqrt{\lambda}}^{-1}\underline{K})^\circ, s{\cal Q}_{{\cal E}}\bigr) N(s{\cal Q}_{{\cal E}}, tB_2^n)
\end{multline*}
for all $t \gs (\alpha^{-1}C_0)^{\frac{2+3\alpha}{2\alpha}}\,C_\alpha^{2/\alpha}$ (and the auxiliary parameter $s$ chosen as before).

{\bf Step 4.} Note that
\begin{multline*}
\max\bigl\{N\bigl(\Delta_{\!\sqrt{\lambda}}^{-1}K, t B_2^n\bigr),N\bigl((\Delta_{\!\sqrt{\lambda}}^{-1}K)^\circ, t B_2^n\bigr), N\bigl(B_2^n, t\Delta_{\!\sqrt{\lambda}}^{-1}K\bigr), N\bigl(B_2^n, t(\Delta_{\!\sqrt{\lambda}}^{-1}K)^\circ\bigr)\bigr\}
\\
\ls \max\bigl\{N\bigl(\Delta_{\!\sqrt{\lambda}}^{-1}\overline{K}, tB_2^n\bigr), N\bigl((\Delta_{\!\sqrt{\lambda}}^{-1}\underline{K})^\circ, tB_2^n\bigr), N\bigl(B_2^n, t\Delta_{\!\sqrt{\lambda}}^{-1}\underline{K}\bigr), N\bigl(B_2^n, t(\Delta_{\!\sqrt{\lambda}}^{-1}\overline{K})^\circ\bigr)\bigr\}.
\end{multline*}

Given our choice of $\alpha$, we have that $\frac{2\alpha}{4+3\alpha} = \beta$. We set 
\begin{equation*}
D_\beta:= 2(\alpha^{-1}\tilde{C}_0)^{\frac{2+3\alpha}{4+3\alpha}}\, C_\alpha^{\frac{4}{4+3\alpha}},
\end{equation*}
and we can then write
\begin{equation}\label{eq:covering-num-total}
\max\bigl\{N\bigl(\Delta_{\!\sqrt{\lambda}}^{-1}K, t B_2^n\bigr),N\bigl((\Delta_{\!\sqrt{\lambda}}^{-1}K)^\circ, t B_2^n\bigr), N\bigl(B_2^n, t\Delta_{\!\sqrt{\lambda}}^{-1}K\bigr), N\bigl(B_2^n, t(\Delta_{\!\sqrt{\lambda}}^{-1}K)^\circ\bigr)\bigr\} \ls \exp\bigl(D_\beta\, n/t^\beta\bigr)
\end{equation}
for all $t\gs D_\beta^{1/\beta}$. Having started with the assumptions that $\underline{K}$ is in $\alpha$-regular $M$-position, while $\Delta_\lambda B_2^n$ is an $\alpha$-regular $M$-ellipsoid for $\overline{K}$, we conclude that the linear image $\Delta_{\!\sqrt{\lambda}}^{-1}K$ of $K$ is in $\beta$-regular $M$-position.

\bigskip

This completes the proof of the theorem in the case that $K$ has Santaló point at the origin.
Note however that \eqref{eq:covering-num-total} is symmetric in $K$ and $K^\circ$ (in the sense that $N\bigl((\Delta_{\!\sqrt{\lambda}}^{-1}K)^\circ, t B_2^n\bigr) = N\bigl(\Delta_{\!\sqrt{\lambda}}(K^\circ), t B_2^n\bigr)$ and $N\bigl(\Delta_{\!\sqrt{\lambda}}^{-1}K, t B_2^n\bigr) = N\bigl((\Delta_{\!\sqrt{\lambda}}(K)^\circ)^\circ, t B_2^n\bigr)$ for instance),
and recall that $K^\circ$ runs over all centred convex bodies as $K$ runs over all convex bodies with 
Santaló point at the origin. Thus we are done in the case of a centred convex body too. 
\qed

\medskip
\bigskip

\noindent{\it Proof of Proposition \ref{prop:Gelfand-numbers}.} Again we start with the assumptions that $K$ has Santaló point at the origin, and that $\underline{K}$ is in $\alpha$-regular $M$-position, where $\alpha=\frac{4\beta}{2-3\beta}$. We now understand this to mean that Pisier's stronger statement holds true: for all $1\ls l\ls n$,
\begin{equation*}
\max\bigl\{c_l\bigl(\underline{K},B_2^n\bigr),\,c_l\bigl((\underline{K})^\circ, B_2^n\bigr)\bigr\}\ls C_0C_\alpha^{1/\alpha}\left(\frac{n}{l}\right)^{1/\alpha}.
\end{equation*}
Recall that this also implies the bounds for the covering numbers of $\underline{K}$ and $(\underline{K})^\circ$ by dilates of $B_2^n$ that we stated before. 

We also find an $\alpha$-regular $M$-ellipsoid ${\cal E}$ for $\overline{K}$ ``in the strong sense'', that is, satisfying regularity for the Gelfand numbers:
\begin{equation*}
\max\bigl\{c_l\bigl(\overline{K},{\cal E}\bigr),\,c_l\bigl((\overline{K})^\circ, {\cal E}^\circ\bigr)\bigr\}\ls C_0C_\alpha^{1/\alpha}\left(\frac{n}{l}\right)^{1/\alpha}
\end{equation*}
(we can find this by first considering an $\alpha$-regular $M$-position of $\overline{K}$ ``in the strong sense'', which we know exists, and then by transforming back to the current position of $\overline{K}$; note that the bounds will remain unchanged since $c_l(TK_1, TL_1)=c_l(K_1,L_1)$ for any linear transformation $T$). Again the bounds for the covering numbers involving $\overline{K}$ that we stated before also hold, so we can conclude that $v_l({\cal E})\ls \tilde{C}\,C_\alpha^{2/\alpha}\bigl(n/l\bigr)^{(2+3\alpha)/\alpha}$ for all $1\ls l\ls n$.
%\begin{equation*}
%v_l({\cal E}) 
%%= \max_{F\in G_{n,l}}{\rm vrad}({\rm Proj}_F({\cal E})) 
%\ls \tilde{C}\,C_\alpha^{2/\alpha} \left(\frac{n}{l}\right)^{(2+3\alpha)/\alpha}
%\end{equation*}
%for all $1\ls l\ls n$.

\medskip

Finally, for convenience again, we can also assume that ${\cal E} = \Delta_\lambda B_2^n$ with $\Delta_\lambda$ the diagonal matrix from the proof of Theorem \ref{main-result}. As before, we set ${\cal Q}_{{\cal E}}= \Delta_{\!\sqrt{\lambda}}B_2^n = \Delta_{\!\sqrt{\lambda}}^{-1}{\cal E}$, and we have
\begin{equation}\label{eqp-main:prop:Gelfand-numbers}
v_l({\cal Q}_{\cal E}) = \sqrt{v_l({\cal E})} \ls \tilde{C}\,C_\alpha^{1/\alpha} \left(\frac{n}{l}\right)^{\frac{2+3\alpha}{2\alpha}}.
\end{equation}
Moreover,
\begin{equation*}
\max\bigl\{ c_l\bigl(\Delta_{\!\sqrt{\lambda}}^{-1}\overline{K},{\cal Q}_{\cal E}\bigr), c_l\bigl((\Delta_{\!\sqrt{\lambda}}^{-1}\overline{K})^\circ, {\cal Q}_{{\cal E}}^\circ\bigr),\ c_l\bigl(\Delta_{\!\sqrt{\lambda}}^{-1}\underline{K},{\cal Q}_{\cal E}^\circ\bigr), c_l\bigl((\Delta_{\!\sqrt{\lambda}}^{-1}\underline{K})^\circ, {\cal Q}_{{\cal E}}\bigr)\bigr\}\ls C_0C_\alpha^{1/\alpha}\left(\frac{n}{l}\right)^{1/\alpha}
\end{equation*}
for all $1\ls l\ls n$.

We start with $l=1$: by the bound for $c_1\bigl(\Delta_{\!\sqrt{\lambda}}^{-1}\overline{K},{\cal Q}_{\cal E}\bigr)$ and by \eqref{eqp-main:prop:Gelfand-numbers}, we can write
\begin{align*}
\Delta_{\!\sqrt{\lambda}}^{-1}K \subseteq \Delta_{\!\sqrt{\lambda}}^{-1}\overline{K}&\subseteq C_0(C_\alpha n)^{1/\alpha} {\cal Q}_{\cal E}
\\
& \subseteq C_0(C_\alpha n)^{1/\alpha} \cdot \sqrt{\lambda_1}B_2^n \subseteq \tilde{C}_0(C_\alpha n)^{1/\alpha} \cdot C_\alpha^{1/\alpha} n^{\frac{2+3\alpha}{2\alpha}}\, B_2^n\subseteq \tilde{\tilde{C}}_0 D_\beta^{1/\beta} n^{1/\beta}\,B_2^n.
\end{align*}

Analogously, we have
\begin{equation*}
\Delta_{\!\sqrt{\lambda}}K^\circ = (\Delta_{\!\sqrt{\lambda}}^{-1}K)^\circ \subseteq (\Delta_{\!\sqrt{\lambda}}^{-1}\underline{K})^\circ \subseteq C_0(C_\alpha n)^{1/\alpha} {\cal Q}_{\cal E} \subseteq \tilde{C}_0 D_\beta^{1/\beta} n^{1/\beta}\,B_2^n.
\end{equation*}

Consider now an even $l\in \{2,\ldots,n\}$. By the bound for $c_{l/2}\bigl(\Delta_{\!\sqrt{\lambda}}^{-1}\overline{K},{\cal Q}_{\cal E}\bigr)$, we can find a subspace $F_0\in G_{n,n-\frac{l}{2}+1}$ such that
\begin{equation*}
(\Delta_{\!\sqrt{\lambda}}^{-1}\overline{K}) \cap F_0\subseteq C_0C_\alpha^{1/\alpha}\left(\frac{2n}{l}\right)^{1/\alpha}\bigl({\cal Q}_{\cal E}\cap F_0\bigr).
\end{equation*}
Let $\mu_1\gs \mu_2\gs\cdots\gs\mu_{n-\frac{l}{2}+1}$ be the lengths of the semiaxes of the ellipsoid ${\cal Q}_{\cal E}\cap F_0$, and let $F_1\ls F_0$ be the subspace spanned by the $n-l+1$ shortest of those (thus the orthogonal complement $F_1^\perp \cap F_0$ of $F_1$ within $F_0$ is the subspace spanned by the $l/2$ longest semiaxes of ${\cal Q}_{\cal E}\cap F_0$). Then ${\cal Q}_{\cal E}\cap F_1\subseteq \mu_{\frac{l}{2}+1}B_{F_1}$ and
\begin{equation*}
\mu_{\frac{l}{2}+1} \ls \biggl(\prod_{i=1}^{l/2}\mu_i\biggr)^{2/l} = {\rm vrad}({\cal Q}_{\cal E}\cap F_1^\perp\cap F_0) \ls w_{l/2}({\cal Q}_{\cal E}) = v_{l/2}({\cal Q}_{\cal E}).
\end{equation*}
Thus, by \eqref{eqp-main:prop:Gelfand-numbers} we obtain
\begin{multline*}
(\Delta_{\!\sqrt{\lambda}}^{-1} K) \cap F_1 \subseteq (\Delta_{\!\sqrt{\lambda}}^{-1}\overline{K}) \cap F_1 \subseteq 
C_0C_\alpha^{1/\alpha}\left(\frac{2n}{l}\right)^{1/\alpha}\bigl({\cal Q}_{\cal E}\cap F_1\bigr)
\subseteq 
C_0C_\alpha^{1/\alpha}\left(\frac{2n}{l}\right)^{1/\alpha}\!\!\cdot\, v_{l/2}({\cal Q}_{\cal E})\, B_{F_1} 
\\
\subseteq 
\tilde{C}_0C_\alpha^{1/\alpha}\left(\frac{2n}{l}\right)^{1/\alpha}\!\!\cdot C_\alpha^{1/\alpha}\left(\frac{2n}{l}\right)^{\frac{2+3\alpha}{2\alpha}}\,B_{F_1}\subseteq \tilde{\tilde{C}}_0 D_\beta^{1/\beta}\left(\frac{n}{l}\right)^{1/\beta}\,B_{F_1}.
\end{multline*} 
This shows that $c_l\bigl(\Delta_{\!\sqrt{\lambda}}^{-1}K,B_2^n\bigr)\ls \tilde{\tilde{C}}_0 D_\beta^{1/\beta}\left(n/l\right)^{1/\beta}$.

\smallskip

Similarly we bound $c_l\bigl(\Delta_{\!\sqrt{\lambda}}K^\circ,B_2^n\bigr)$ by using the bound for $c_{l/2}\bigl((\Delta_{\!\sqrt{\lambda}}^{-1}\underline{K})^\circ, {\cal Q}_{{\cal E}}\bigr)$ and \eqref{eqp-main:prop:Gelfand-numbers}.

Finally, we observe that the sequences $l\in \{1,2,\ldots,n\}\mapsto c_l\bigl(K_1, B_2^n)$, where $K_1=\Delta_{\!\sqrt{\lambda}}^{-1}K$ or $=\Delta_{\!\sqrt{\lambda}}K^\circ$, are decreasing, so by slightly adjusting the absolute constant $C_0$, we can also bound the remaining terms of those sequences. 

As before, if $K$ has barycentre at the origin instead, then we work along the same lines, but with $K^\circ$ taking on the role of $K$. This completes the proof. \qed

\bigskip

Recall that to prove Theorem \ref{main-result}, we relied on \eqref{eqp0:thm:main-result}, which gives a comparison of the volumes of corresponding projections of $\overline{K}$ and of $\underline{K}$ when $s(K)=0$ (with estimates which depend on the dimensions of the projections in a regular way). Using Theorem \ref{main-result} now, we can establish a similar comparison for volumes of such projections when $K$ is centred.

\begin{corollary}\label{cor:main-result}
Let $K$ be a centred convex body in ${\mathbb R}^n$, and let $1\ls l\ls n$. For every $F\in G_{n,l}$ we have
\begin{equation*}
|{\rm Proj}_F(\overline{K})|^{1/l}\ls C \,\Bigl(\frac{n}{l}\Bigr)^5\bigl(\log(en/l)\bigr)^2\,|{\rm Proj}_F(\underline{K})|^{1/l}. 
\end{equation*}
\end{corollary}

We defer the proof of the corollary to the next section.

\section{Remarks and applications}

1. The proof of Theorem \ref{main-result} shows that, if $K\subset {\mathbb R}^n$ has Santaló point at the origin,
then a $\beta$-regular position of $K$ can be found if we apply to $K$ the linear transformation $T\in {\rm GL}(n)$
that takes $\underline{K}$ to $\alpha$-regular $M$-position 
%in the sense of Theorem \ref{thm:Pisier-regular-ellipsoids},
and then adjust a little so that $\widetilde{\,\,\overline{K}\,\,}$ and $\bigl(\widetilde{\,\,\underline{K}\,\,}\bigr)^\circ$ are covered in a regular way by the same ellipsoid. Similarly, if $K$ has barycentre at the origin, then we start by placing $\overline{K}$ in $\alpha$-regular $M$-position, where $\alpha = \frac{4\beta}{2-3\beta}$.

It seems very likely that the regularity of the covering numbers 
%for an arbitrary non-symmetric convex body 
given by this method
is not optimal, and one could even conjecture that there exists
a $1$-regular $M$-position for any not-necessarily symmetric convex body.
That said, if we were to keep the core principle of this method, which is to rely on an inequality for volumes of projections in the spirit of \eqref{eqp0:thm:main-result}, then most probably we wouldn't be able to achieve this. This is because the example of the simplex shows that the exponent of $\frac{n}{l}$ in any inequality of similar form to \eqref{eqp0:thm:main-result} should be at least $\frac{1}{2}$ (which plugged into the rest of our argument would give at best $\frac{6}{11}$-regularity).

Indeed, assume $l\ls \frac{n}{2}$. As we have already seen, the regular simplex $\hat{S}_n$ with edge-length $\frac{n+1}{\sqrt{2}}$ has $l$-dimensional sections through its barycentre with volume
\begin{equation*}
\Bigl(\frac{n+1}{2}\Bigr)^l\frac{\sqrt{n+1\,}}{l!\,\sqrt{n+1-l\,}}
\end{equation*}
Again we recall that this simplex can be embedded in ${\mathbb R}^{n+1}$ so that its barycentre is at the origin, and so that the hyperplane containing it is the hyperplane $H_{\bm 1}$ orthogonal to the vector ${\bm 1} = e_1+e_2+\ldots+e_{n+1}$. Then 
\begin{equation*}
\hat{S}_n\cap(-\hat{S}_n) = \tfrac{1}{2}B_\infty^{n+1} \cap H_{\bm 1}
%{\rm Cube}_{n+1}\cap H_{\bm 1}
\end{equation*}
where $\tfrac{1}{2}B_\infty^{n+1} = \bigl[-\frac{1}{2},\frac{1}{2}\bigr]^{n+1}$ is the (origin-symmetric) cube of volume 1 in ${\mathbb R}^{n+1}$ (see e.g. \cite{Taschuk-thesis} where the idea for this normalisation is taken from). 

Consider now any $l$-dimensional subspace $F$ of the hyperplane $H_{\bm 1}$; then 
\begin{equation*}
\bigl(\hat{S}_n\cap(-\hat{S}_n)\bigr)\cap F = \bigl(\tfrac{1}{2}B_\infty^{n+1} \cap H_{\bm 1}\bigr)\cap F = \tfrac{1}{2}B_\infty^{n+1} \cap F.
\end{equation*}
It remains to recall a result by K. Ball \cite{Ball-1989} (optimal in many cases) which gives that, for any such subspace $F$,
\begin{equation*}
|\tfrac{1}{2}B_\infty^{n+1} \cap F|_l\ls \biggl(\sqrt{\frac{n+1}{l}}\biggr)^l.
\end{equation*}
It follows that we cannot have
\begin{equation*}
\big|\bigl(\hat{S}_n-\hat{S}_n\bigr)\cap F\big|^{1/l}\ls C\Bigl(\frac{n}{l}\Bigr)^\gamma\,\big|\bigl(\hat{S}_n\cap(-\hat{S}_n)\bigr)\cap F\big|^{1/l}
\end{equation*}
with $\gamma < \frac{1}{2}$ (note also that by polarity and the Blaschke-Santaló and Bourgain-Milman inequalities, this last inequality is equivalent to an inequality of the form of \eqref{eqp0:thm:main-result}).

\medskip

2. Our first remark notwithstanding, in the case of the simplex there does exist an $\alpha$-regular $M$-position for any $\alpha\in (0,2)$. Still it should be said that it also seems likely this would not be true in full generality. 

To verify this, we now consider the regular simplex $S_n$ of edge-length $\sqrt{2}$ and barycentre at the origin, embedded in ${\mathbb R}^{n+1}$:
\begin{equation*}
S_n = {\rm conv}\Bigl\{e_i-\frac{1}{n+1}\bm{1}:1\ls i\ls n+1\Bigr\}.
\end{equation*}
%where $\{e_i: 1\ls i\ls n+1\}$ is the standard orthonormal basis in ${\mathbb R}^{n+1}$
%and ${\bm 1}\in {\mathbb R}^{n+1}$ is the vector all of whose coordinates equal 1. 
Again recall that $\overline{S_n} = {\rm conv}\{S_n, -S_n\}$ can be viewed as a projection of the $\ell_1$-unit ball $B_1^{n+1}$:
\begin{equation*} 
{\rm conv}\{S_n, -S_n\} = {\rm Proj}_{H_{\bm 1}}\bigl(B_1^{n+1}\bigr).
\end{equation*}
%where $E={\bm 1}^\perp$ is the hyperplane of vectors in ${\mathbb R}^{n+1}$ whose coordinates add to 0. 

Recall also that $S_n$ satisfies $S_n^\circ = -(n+1)S_n$
(where the polar set is considered within $H_{\bm 1}$), therefore
\begin{align*}
\underline{S_n} = S_n\cap (-S_n) = \frac{1}{n+1}\bigl(S_n^\circ\cap (-S_n^\circ)\bigr)
&= \frac{1}{n+1}\bigl({\rm conv}\{S_n, -S_n\}\bigr)^\circ 
\\&= \frac{1}{n+1}\bigl({\rm Proj}_{H_{\bm 1}}\bigl(B_1^{n+1}\bigr)\bigr)^\circ
 = \frac{1}{n+1} \bigl(B_\infty^{n+1}\cap H_{\bm 1}\bigr).
\end{align*} 

Schütt \cite{Schutt-1984} has shown that
\begin{equation*}
e_{k+1}\Bigl(B_1^{n+1}, \frac{1}{\sqrt{n}}B_2^{n+1}\Bigr) \simeq \left\{\begin{array}{cl}
\sqrt{n} &\  \hbox{if}\ \,k\ls \log(n+1)
\smallskip \\ \sqrt{\dfrac{n}{k}}\,\sqrt{\log(en/k)} &\  \hbox{if}\ \,\log(n+1)\ls k\ls n+1
\smallskip \\ 2^{-k/n}n &\  \hbox{if} \ \,k\gs n+1
\end{array}\right..
\end{equation*}
We combine this with the fact that, for every $p > \frac{1}{2}$ and any $1\ls k\ls n+1$,
\begin{gather*} 
\exp\left(\Bigl(p-\frac{1}{2}\Bigr)\log(en/k)\right) \gs 
%\exp\left(\frac{p-1/2}{2}\,\log(en/k)\right) =
\sqrt{\exp\left(\Bigl(p-\frac{1}{2}\Bigr)\log(en/k)\right)} \gs \sqrt{\Bigl(p-\frac{1}{2}\Bigr)\log(en/k)}
\\
\Longrightarrow\quad 
e_{k+1}\Bigl(B_1^{n+1}, \frac{1}{\sqrt{n}}B_2^{n+1}\Bigr) \lesssim \sqrt{\dfrac{en}{k}}\,\sqrt{\log(en/k)} \ls \frac{1}{\sqrt{p-\frac{1}{2}}}\,\left(\frac{en}{k}\right)^p \qquad
\end{gather*}
to obtain that
%for $t\gs \bigl(p-\tfrac{1}{2}\bigr)^{-1/2}$,
\begin{equation*}
%e_{k+1}\Bigl(B_1^{n+1}, \frac{1}{\sqrt{n}}B_2^{n+1}\Bigr) \lesssim \frac{e^p}{\sqrt{p-\frac{1}{2}}}\,\left(\frac{n}{k}\right)^p
%\quad \Longrightarrow \quad 
N\Bigl(B_1^{n+1}, \frac{t}{\sqrt{n}}B_2^{n+1}\Bigr) \ls \exp\Bigl(C\bigl(p-\tfrac{1}{2}\bigr)^{-\frac{1}{2p}}\,n/t^{1/p}\Bigr)
\end{equation*}
for any $t\gs e^p\bigl(p-\tfrac{1}{2}\bigr)^{-1/2}$, with $C$ an absolute constant.

At the same time, by the dual Sudakov inequality (see e.g. \cite[Theorem 4.2.2]{AGM-book})
\begin{equation*}
N\Bigl(\frac{1}{\sqrt{n}}B_2^{n+1}, tB_1^{n+1}\Bigr) \ls \exp (Cn/t^2)
\end{equation*}
for every $t>0$. 

As a consequence,
\begin{equation*} 
\max\Bigl\{N\Bigl({\rm Proj}_{H_{\bm 1}}\bigl(B_1^{n+1}\bigr), \frac{t}{\sqrt{n}}B_2^n\Bigr), 
N\Bigl(\frac{1}{\sqrt{n}}B_2^n, t{\rm Proj}_{H_{\bm 1}}\bigl(B_1^{n+1}\bigr)\Bigr)\Bigr\}
\ls \exp\Bigl(C^\prime\bigl(p-\tfrac{1}{2}\bigr)^{-\frac{1}{2p}}\,n/t^{1/p}\Bigr)
\end{equation*}
for every $t\gs e^p\bigl(p-\tfrac{1}{2}\bigr)^{-1/2}$,
and by the duality of covering numbers \cite{AMS-2004} the same bound
(with slightly different absolute constants) is valid for
\begin{equation*}
\max\Bigl\{N\Bigl(\bigl({\rm Proj}_{H_{\bm 1}}\bigl(B_1^{n+1}\bigr)\bigr)^\circ, t\sqrt{n}B_2^n\Bigr), 
N\Bigl(\sqrt{n}B_2^n, t\bigl({\rm Proj}_{H_{\bm 1}}\bigl(B_1^{n+1}\bigr)\bigr)^\circ\Bigr)\Bigr\}.
\end{equation*}

\smallskip

We thus see that both $\sqrt{n}\,\overline{S_n}$ and $\sqrt{n}\,\underline{S_n}$ are in $\alpha$-regular $M$-position
for any $\alpha\in (0,2)$ with a constant of the same form as in Theorem \ref{thm:Pisier-regular-ellipsoids}
(in fact, here it suffices to know that one of them is in regular position since $\bigl(\,\overline{S_n}\,\bigr)^\circ = (n+1)\underline{S_n}$).
This now obviously implies that $\sqrt{n}\,S_n$ has the same property as well.

\medskip

3. It should be clear from Subsection \ref{subsec:KM-approach} that 
Klartag and Milman's argument could also be modified to work in the general case
if there were a version of Fradelizi's result: 
\begin{equation} \label{eq:modified-Fradelizi-result1}
\min\limits_{F\in G_{n,l}} \frac{|K\cap F|}{\max\limits_{y\in {\mathbb R}^n}|K\cap (y+F)|}\gs \Bigl(\frac{l+1}{n+1}\Bigr)^l 
\end{equation}
also available when $K$ has Santaló point at the origin.

Such a version, albeit with an estimate which, we should probably expect, is rather crude, follows now by Theorem \ref{main-result}. 
%and Proposition \ref{prop:BS-ineq-proj-sec}. 
Indeed, let $K$ be a body with Santaló point at the origin and $1\ls l < n$.
Since the left-hand side of \eqref{eq:modified-Fradelizi-result1}
is invariant under linear transformations, we may assume that $K$ is in $\beta_l$-regular position with
\begin{equation*} \beta_l = \frac{2}{5} - \frac{1}{5\log(en/l)}. \end{equation*}
Then, for every $F\in G_{n,l}$,
\begin{equation*} 
\max\bigl\{|{\rm Proj}_F(K)|^{1/l}, |{\rm Proj}_F(K^\circ)|^{1/l}\bigr\} \ls e D_{\beta_l}^{1/\beta_l}\Bigl(\frac{n}{l}\Bigr)^{1/\beta_l} |B_2^l|^{1/l}.
\end{equation*}
Combining this with the Bourgain-Milman inequality, we see that
\begin{align*}
|K\cap F|^{1/l} & \gs c_0\frac{|B_2^l|^{2/l}}{|{\rm Proj}_F(K^\circ)|^{1/l}}  \gs \frac{c_0}{eD_{\beta_l}^{1/\beta_l}}
\Bigl(\frac{l}{n}\Bigr)^{1/\beta_l} |B_2^l|^{1/l}
\\
& \gs \frac{c_0}{e^2D_{\beta_l}^{2/\beta_l}} \Bigl(\frac{l}{n}\Bigr)^{2/\beta_l}|{\rm Proj}_F(K)|^{1/l}
\gs \frac{c_0}{e^2D_{\beta_l}^{2/\beta_l}} \Bigl(\frac{l}{n}\Bigr)^{2/\beta_l} \max\limits_{y\in {\mathbb R}^n}|K\cap (y+F)|^{1/l}.
\end{align*}
Noting that $D_{\beta_l}^{2/\beta_l} \lesssim \bigl(\log(en/l)\bigr)^2$, we conclude that
\begin{equation*}
\min\limits_{F\in G_{n,l}} \frac{|K\cap F|}{\max\limits_{y\in {\mathbb R}^n}|K\cap (y+F)|}
\gs \frac{c^l}{\bigl(\log(en/l)\bigr)^{2l}}\Bigl(\frac{l}{n}\Bigr)^{5l} 
%\gs C^{-\max\{l,n-l\}}
\end{equation*}
for some absolute constant $c>0$.

\subsection{Comparison of volumes of projections of $\overline{K}$ and $\underline{K}$ when $K$ is centred}

%With an analogue of Fradelizi's result for convex bodies with Santaló point at the origin, we can now compare $v_l^-(\underline{K})$ with $v_l^-(\overline{K})\simeq v_l^-(K)$. Repeating the steps from Subsection \ref{subsec:KM-approach}, except for replacing Fradelizi's result with the analogue we gave above, we get for every $F\in G_{n,l}$ that
%\begin{equation}\label{eq:proj2sec-Santalo}
%{\rm vrad}({\rm Proj}_F(\underline{K}))
%\gs \frac{c}{\log^2(en/l)}\left(\frac{l}{n}\right)^7{\rm vrad}({\rm Proj}_F(K))
%\gs \frac{c}{\log^2(en/l)}\left(\frac{l}{n}\right)^7v_l^-(K).
%\end{equation}

By repeating the steps from Subsection \ref{subsec:KM-approach}, with the only modification being that we replace Fradelizi's result with the analogue we obtained above, we could only get $|({\rm Proj}_F(\overline{K}))|^{1/l} \ls C(n/l)^7{\log^2(en/l)} |({\rm Proj}_F(\underline{K}))|^{1/l}$ for $F\in G_{n,l}$. Instead, we can get a better estimate if we apply the 1st step of the proof of Theorem \ref{main-result} carefully.

\medskip

\noindent
{\it Proof of Corollary \ref{cor:main-result}.} Fix $1\ls l < n$, and set $\alpha_l = 2-\frac{1}{\log(en/l)}$. Since 
\begin{equation*}
\sup_{F\in G_{n,l}} \frac{|{\rm Proj}_F(\overline{K})|^{1/l}}{|{\rm Proj}_F(\underline{K})|^{1/l}} 
\end{equation*}
is linearly invariant, we can choose a suitable position for the centred convex body $K$. As seen in the proof of Theorem \ref{main-result}, if $\overline{K}$ is in $\alpha_l$-regular $M$ position, and ${\cal D}$ is an $\alpha_l$-regular $M$-ellipsoid for $\underline{K}$, then 
\begin{equation*}
w_l({\cal D}^\circ) = v_l({\cal D}^\circ) \ls C\,C_{\alpha_l}^{2/\alpha_l} \left(\frac{n}{l}\right)^{(2+3\alpha_l)/\alpha_l}.
\end{equation*}
We can then write
\begin{align*}
|{\rm Proj}_F(\overline{K})|^{1/l} &\ls C\,C_{\alpha_l}^{1/\alpha_l}\left(\frac{n}{l}\right)^{1/\alpha_l}\,|B_2^l|^{1/l}
\\
&\ls C\,\frac{1}{v_l^-({\cal D})}\,C_{\alpha_l}^{1/\alpha_l}\left(\frac{n}{l}\right)^{1/\alpha_l}\,|{\rm Proj}_F({\cal D})|^{1/l}
\\
&= C\,w_l({\cal D}^\circ)\,C_{\alpha_l}^{1/\alpha_l}\left(\frac{n}{l}\right)^{1/\alpha_l}\,|{\rm Proj}_F({\cal D})|^{1/l}
\\
&\ls C^\prime\,C_{\alpha_l}^{3/\alpha_l} \left(\frac{n}{l}\right)^{(3+3\alpha_l)/\alpha_l} \,|{\rm Proj}_F({\cal D})|^{1/l}
\\
&\ls C^{\prime\prime}\,C_{\alpha_l}^{4/\alpha_l} \left(\frac{n}{l}\right)^{(4+3\alpha_l)/\alpha_l}\,|{\rm Proj}_F(\underline{K})|^{1/l}.
\end{align*}
The conclusion follows. \qed

\begin{remark}
We could also use Theorem \ref{main-result} as a black box (that is, rely only on its statement), but we gave the above proof first, which indicates that the almost $\frac{2}{5}$-regularity that we got is already hinted at in applications of the first step of the proof of the theorem. If instead we start with a linear image $K_l$ of the centred convex body $K$ which satisfies
\begin{equation*}
\max\left\{N(K_l,tB_2^n),\,N(K_l^\circ, tB_2^n)\right\}\ls \exp\left(D_{\beta_l}n/t^{\beta_l}\right)
\end{equation*}
for all $t\gs D_{\beta_l}^{1/\beta_l}$ with $\beta_l = \frac{2}{5} - \frac{1}{5\log(en/l)}$, then we can observe the following: since $K_l$ is convex and contains the origin in its interior, we have $\overline{K_l}={\rm conv}(K_l,-K_l)=\{\lambda x+(1-\lambda)y:x,-y\in K_l,\lambda\in [0,1]\} \subset K_l-K_l$, and hence
\begin{equation*}
N(\overline{K_l},2tB_2^n)\ls N(K_l-K_l,2tB_2^n)\ls \bigl(N(K_l, tB_2^n)\bigr)^2 \ls \exp\left(2D_{\beta_l}n/t^{\beta_l}\right)
\end{equation*}
for $t$ in the same range. In the same way, we get
\begin{equation*}
N(\overline{K_l^\circ},2tB_2^n)\ls  \bigl(N(K_l^\circ, tB_2^n)\bigr)^2 \ls \exp\left(2D_{\beta_l}n/t^{\beta_l}\right), 
\end{equation*}
which by the duality of covering numbers (see \cite{AMS-2004}) gives
\begin{equation*}
N(B_2^n, \gamma_1 t\,\underline{K_l}) = N\bigl(B_2^n, \gamma_1 t \bigl(\overline{K_l^\circ}\bigr)^\circ\bigr) \ls \exp\left(\gamma_2 D_{\beta_l}n/t^{\beta_l}\right)
\end{equation*}
for some absolute constants $\gamma_1, \gamma_2$, and for $t$ in the same range. We can now write, for any $F\in G_{n,l}$,
\begin{equation*}
|{\rm Proj}_F(\overline{K_l})|^{1/l} \ls C\,D_{\beta_l}^{1/\beta_l}\left(\frac{n}{l}\right)^{1/\beta_l}\,|B_2^l|^{1/l}\ls C^\prime_{\gamma_1,\gamma_2} \,D_{\beta_l}^{2/\beta_l}\left(\frac{n}{l}\right)^{2/\beta_l}\,|{\rm Proj}_F(\underline{K_l})|^{1/l},
\end{equation*}
which by our choice of $\beta_l$ leads to the same estimates as above.

With this we can also conclude the necessity and sufficiency of an inequality such as \eqref{eqp0:thm:main-result} in order to have regular covering of a non-symmetric convex body $K$ (with $0\in {\rm int}(K)$): by Klartag and Milman's suggested approach, such an inequality (for the specific body $K$) leads to regular $M$-ellipsoids for $K$, while if we know that $K$ has $\beta$-regularity, or almost $\beta$-regularity, then an inequality such as \eqref{eqp0:thm:main-result} follows (with the exponent of the factor $n/l$ being $2/\beta$, or $2/\beta+o(1)$). Of course we don't have precise equivalence, since the estimates become worse as we pass from one result to the other.
\end{remark}

\medskip

As an aside, we also give two more estimates for $\sup_{F\in G_{n,l}} \frac{|{\rm Proj}_F(\overline{K})|^{1/l}}{|{\rm Proj}_F(\underline{K})|^{1/l}}$ when $K$ is centred, which could be useful in specific cases.

\smallskip

a. We have that
\begin{equation*}
\frac{|B_2^l|^{2/l}}{|{\rm Proj}_F(\underline{K})|^{1/l}}\ls C|\overline{K^\circ}\cap F|^{1/l}\ls C^\prime\left(\frac{n}{l}\right)\max_{x\in {\mathbb R}^n}|K^\circ\cap(x+F)|^{1/l}\ls 
C^{\prime\prime}\left(\frac{n}{l}\right)^2\big|\bigl(K^\circ-b(K^\circ)\bigr)\cap F\big|^{1/l}
\end{equation*}
by the Bourgain-Milman inequality, Rudelson's and Fradelizi's results. 

We now observe that
\begin{equation*}
K^\circ-b(K^\circ)\subset K^{\circ}+\|-b(K^\circ)\|_{K^\circ}K^\circ = [1+\|-b(K^\circ)\|_{K^\circ}]K^\circ = [1+h_K(-b(K^\circ))]K^\circ.
\end{equation*}
Thus
\begin{align*}
\left(\frac{n}{l}\right)^2\big|\bigl(K^\circ-b(K^\circ)\bigr)\cap F\big|^{1/l} &\ls [1+h_K(-b(K^\circ))]\left(\frac{n}{l}\right)^2|K^\circ\cap F|^{1/l} 
\\
&\ls C[1+h_K(-b(K^\circ))]\left(\frac{n}{l}\right)^3 \frac{|B_2^l|^{2/l}}{|{\rm Proj}_F(K)|^{1/l}} 
\\[0.3em]
&\ls C^\prime
\,[1+h_K(-b(K^\circ))]\left(\frac{n}{l}\right)^3 \frac{|B_2^l|^{2/l}}{|{\rm Proj}_F(\overline{K})|^{1/l}}
\end{align*}
with the penultimate inequality due to Theorem \ref{thm:BS-ineq-proj-sec} (or the argument by Klartag and Milman in Subsection \ref{subsec:KM-approach}).

We conclude that
\begin{equation}\label{eq:proj2sec-centred-h}
\sup_{F\in G_{n,l}} \frac{|{\rm Proj}_F(\overline{K})|^{1/l}}{|{\rm Proj}_F(\underline{K})|^{1/l}}\ls C_0[1+h_K(-b(K^\circ))]\left(\frac{n}{l}\right)^3
\end{equation}
for some absolute constant $C_0>0$ (note that this statement also recovers the corresponding estimate for convex bodies $K$ with Santaló point at the origin, which was to be expected since we tried to emulate that one).

\medskip

b. Since $\sup_{F\in G_{n,l}} \frac{|{\rm Proj}_F(\overline{K})|^{1/l}}{|{\rm Proj}_F(\underline{K})|^{1/l}}$ is linearly invariant, we may assume that $K$ is in isotropic position. Using both types of Rogers-Shephard inequalities that we stated (inequality \eqref{eq:ineq-sym-conv} and inequality \eqref{eq:RS-Spingarn}), we have that, for every $F\in G_{n,l}$,
\begin{equation*}
|{\rm Proj}_F(\overline{K})|^{1/l}\simeq |{\rm Proj}_F(K)|^{1/l}\ls C\frac{n}{l}\,\frac{1}{|K\cap F^\perp|^{1/l}}.
\end{equation*}
We now combine several of the properties of the isotropic position that we mentioned:
%\begin{equation*}
%|K\cap F^\perp|^{1/l}\simeq \frac{L_{\pi_F({\bm 1}_K)}}{L_K}
%\end{equation*}
%and hence 
\begin{equation*}
C\frac{n}{l}\,\frac{1}{|K\cap F^\perp|^{1/l}}\ls C^\prime \frac{n}{l}\,\frac{L_K}{L_{\pi_F({\bm 1}_K)}}\ls C^{\prime\prime}\frac{n}{l}L_K
\end{equation*}
(given also that all isotropic constants are uniformly bounded from below),
while $K\supset c_0L_KB_2^n$ for some absolute constant $c_0$ immediately gives $\underline{K}\supset c_0L_KB_2^n$ too. Thus
\begin{equation*}
|{\rm Proj}_F(\overline{K})|^{1/l}\ls C\frac{n}{l}\,\sqrt{l\,}|L_KB_2^l|^{1/l} \ls C\frac{n}{\sqrt{l}}\,|{\rm Proj}_F(\underline{K})|^{1/l}.
\end{equation*}
Despite its apparent crudeness, this is a better estimate, even than \eqref{eq:proj2sec-centred-h}, for small dimensions $l$ (in particular, for $l\ls n^{4/5}$). Still, it has been unclear to us whether it could make a difference in any part of the arguments of this paper.

\subsection{Application to results around the mean norm of isotropic convex bodies}

Aside from its namesake conjecture, the isotropic position is a very useful normalisation that comes up in other problems too, including ones that are motivated by algorithmic applications. It thus appears natural to also examine how the position behaves with respect to parameters which have been introduced in completely different contexts and are optimised in/characterise different positions of a convex body. In certain cases we can observe a substantially worse behaviour compared to the optimal position for the parameter in question (see e.g. \cite{Markessinis-Valettas-2015} and the results there regarding surface area), and in other cases we have a very comparable behaviour.

Here we focus on the corresponding questions regarding the parameters of mean width and mean norm. 
First, let us recall that, by results of Figiel--Tomczak-Jaegermann \cite{Figiel-Tomczak-1979}, Lewis \cite{Lewis-1979} and Pisier \cite{Pisier-1982}, we have that, if $K\subset {\mathbb R}^n$ is origin-symmetric, then 
%there exists $T\in {\rm GL}(n)$ 
\begin{equation}\label{eq:optimal-means-symmetric}
\inf_{T\in {\rm GL}(n)}M(TK)\cdot M^\ast(TK)= \inf_{T\in {\rm GL}(n)}M(TK)\cdot M((TK)^\circ)\ls C\log(n)
\end{equation}
where $C$ is an absolute constant. Also, that
\begin{equation*}
\frac{1}{M(K)}\ls {\rm vrad}(K)\ls M^\ast(K),
\end{equation*}
where the first inequality follows by Jensen's inequality, while the second one is the classical Urysohn inequality (see e.g. \cite[Corollary 1.4]{Pisier-book} for a proof). Thus, for every symmetric convex body $K$ {\bf of volume 1} in ${\mathbb R}^n$, we can find $T_0\in {\rm GL}(n)$ such that
\begin{equation}\label{eq:separate-optimal-means}
M(T_0K)\ls \tilde{C}\frac{\log(n)}{\sqrt{n}} \quad \hbox{and}\ \ M^\ast((T_0K)^\circ)\ls \tilde{C}\sqrt{n}\log(n).
\end{equation}
In the case of a not-necessarily symmetric convex body $L$ with $0\in {\rm int}(L)$, we can instead use the relations $h_L(\theta) = h_{L-0}(\theta) \ls h_{L-L}(\theta) = \max\{\langle x-z,\theta\rangle: x,z\in L\}\ls h_L(\theta) + h_L(-\theta)$ to obtain
\begin{equation*}
M^\ast(L) \ls M^\ast(L-L)\ls 2M^\ast(L),
\end{equation*}
and then we can suitably transform $L-L$ or $L^\circ-L^\circ$ to upper-bound $M^\ast(L)$ or $M(L)=M^\ast(L^\circ)$ respectively by the same, almost optimal, bounds above (that is, if we also assume that $|L|=1$, and in the latter case if we can also say that ${\rm vrad}((L^\circ-L^\circ)^\circ)\simeq {\rm vrad}(L\cap (-L))$ is $\simeq {\rm vrad}(L)$; by the Rogers-Shephard, Blaschke-Santaló and Bourgain-Milman inequalities, this will in fact hold true if $L$ has barycentre of Santaló point at the origin). 

We should also remark that the best known analogue of \eqref{eq:optimal-means-symmetric} in the not-necessarily symmetric case is by Rudelson \cite{Rudelson-2000b} (in a companion paper to \cite{Rudelson-2000a}), who showed that
\begin{equation*}
%\inf_{\substack{T\in {\rm GL}(n) \\ x\in {\rm int}(L)}}M(T(L-x))\cdot M^\ast(T(L-x)) 
\inf\bigl\{M(T(L-x))\cdot M^\ast(T(L-x)) \,:\, T\in {\rm GL}(n),\,x\in {\rm int}(L)\bigr\}
\ls Cn^{1/3}\log^a(n)
\end{equation*}
for some absolute positive constant $a \ls 9$.

Going back to the question of estimating the mean width and the mean norm in the isotropic position, obviously we are no longer allowed to transform $K$ (or $K-K$ or $K^\circ-K^\circ$) further. In the case of mean width, E. Milman \cite{EMilman-2015} showed that, for every isotropic convex body $K$ in ${\mathbb R}^n$,
\begin{equation*}
M^\ast(Z_n(K))\ls C\sqrt{n}\log^2(n) L_K
\end{equation*}
where $Z_n(K)$ is the $L_n$-centroid body of $K$ (its definition is not relevant here, see e.g. \cite{EMilman-2015} or \cite{BGVV-book} for references). One should then recall that $Z_n(K) \simeq K-K$ (see e.g. \cite[Corollary 3.2.9]{BGVV-book}), whence $M^\ast(Z_n(K)) \simeq M^\ast(K-K)\simeq M^\ast(K)$ (regardless of whether $K$ is symmetric or simply centred). In fact, following the recent breakthroughs leading to logarithmic bounds for the isotropic constant, the bound here is now essentially optimal for all practical purposes. 
%(as some logarithmic factor should be there anyway, as evidenced by the normalised $\ell_1$-unit ball, and we usually don't mind some `extra' such factors in applications).

\medskip

In the case of mean norm however, it is not clear at all whether we should expect the symmetric and the non-symmetric case to exhibit comparable behaviours. For symmetric $n$-dimensional isotropic convex bodies, an approach by Giannopoulos and E. Milman \cite{Giannopoulos-EMilman-2014} leads to the following upper bound:
\begin{equation*}
M(K) \ls C\frac{\log^b(n)}{n^{1/6}L_K}
\end{equation*}
where the exponent $b\ls 0.5$ (note that at the time Giannopoulos and Milman didn't have access to the improved bounds for the isotropic constant that we can now plug into their intermediate results). Their approach involved the following main steps.

\medskip

First, Dudley's entropy estimate: 
\begin{equation}\label{eq:Dudley-estimate}
\sqrt{n}M^\ast(K)\ls C\sum_{l=1}^n\frac{1}{\sqrt{l}}e_l(K,B_2^n)
\end{equation}
(see \cite[Theorem 5.5]{Pisier-book} for this formulation, and also note that, due to standard volumetric estimates in an $n$-dimensional space, summation can be taken up to $l=n$). Observe that this estimate holds true for non-symmetric convex bodies as well: we have already seen that $M^\ast(K)\simeq M^\ast(K-K)$. Moreover, for every even integer $l$ we have
\begin{equation*}
e_{l-1}(K-K,B_2^n) \ls 2e_{l/2}(K,B_2^n)
\end{equation*}
(since $N(K-K, 2t_0B_2^n)\ls N(K, t_0 B_2^n)N(-K,t_0 B_2^n)\ls 2^{l-2}\,$ if $\,t_0=e_{l/2}(K,B_2^n)$).

\medskip

Secondly, the following proposition:
\begin{proposition}{\rm (\cite[Theorem 1]{Giannopoulos-EMilman-2014})}\label{prop1:Gian-Milman}
Let $K$ be an origin-symmetric convex body in ${\mathbb R}^n$. For every $l\in \{1,2,\ldots,n\}$ we have
\begin{align*}
e_l(K,B_2^n) &\ls C\frac{n}{l}\log\left(\frac{en}{l}\right)\sup\bigl\{2^{-\frac{l}{3m}}w_m(K):1\ls m\ls l\bigr\}
\\
&\simeq C\frac{n}{l}\log\left(\frac{en}{l}\right)\sup\Bigl\{2^{-\frac{l}{3m}}\frac{1}{v_m^-(K^\circ)}:1\ls m\ls l\Bigr\}.
\end{align*}
\end{proposition}

Note that the second inequality did not need to be stated in \cite[Theorem 1]{Giannopoulos-EMilman-2014} as it is an immediate consequence of the Blaschke-Santaló and Bourgain-Milman inequalities. However, in the non-symmetric case, as we have already seen, we would sometimes have to pay an extra factor of $\frac{n}{l}$ to pass from one quantity to the other one, so in what follows we will aim to show the one of the two inequalities which is `less expensive' to get (and which fortunately is also the most convenient one for our application).

\medskip

Furthermore, Giannopoulos and Milman obtain similar upper bounds for the Gelfand numbers of a symmetric convex body $K$ (which due to Carl's theorem (see Subsection \ref{sec:prelim}) could replace the bounds for the entropy numbers in the application to the mean norm of $K$, and in the end allowed them to slightly simplify those computations).
\begin{proposition}{\rm (\cite[Theorem 2]{Giannopoulos-EMilman-2014})} \label{prop2:Gian-Milman}
Let $K$ be an origin-symmetric convex body in ${\mathbb R}^n$. For every $l\in \{1,2,\ldots,\lfloor\frac{n}{2}\rfloor\}$ we have
\begin{equation*}
c_{2l-1}(K,B_2^n) \ls C\frac{n}{l}\log\left(\frac{en}{l}\right)w_l(K)\simeq C\frac{n}{l}\log\left(\frac{en}{l}\right)\frac{1}{v_l^-(K^\circ)}.
\end{equation*}
\end{proposition}
%
%Again, the second inequality did not need to be stated in their paper, but it is the one we will try to adapt.

We can now give some analogous versions of these two propositions for not-necessarily symmetric convex bodies. 

\begin{proposition}\label{prop:analogue-Gian-Milman1}
Let $K$ be a convex body in ${\mathbb R}^n$, and let $l\in\{1,2,\ldots,n\}$. If $K$ has barycentre at the origin, then
\begin{equation}\label{eq1:prop:analogue-Gian-Milman1}
e_l(K,B_2^n)\ls C\frac{n}{l}\log\left(\frac{en}{l}\right)\sup\Bigl\{2^{-\frac{l}{3m}}\left(\frac{n}{m}\right)^2 w_m(K):1\ls m\ls l\Bigr\}
\end{equation}
and
\begin{equation*}
e_l(K,B_2^n)\ls C\frac{n}{l}\log\left(\frac{en}{l}\right)\sup\Bigl\{2^{-\frac{l}{3m}}\left(\frac{n}{m}\right)^3\frac{1}{v_m^-(K^\circ)}:1\ls m\ls l\Bigr\}.
\end{equation*}

On the other hand, if $K$ has either barycentre or Santaló point at the origin, then
\begin{equation}\label{eq2:prop:analogue-Gian-Milman1}
e_l(K,B_2^n)\ls C\left(\frac{n}{l}\right)^5\log^2\left(\frac{en}{l}\right)\sup\Bigl\{2^{-\frac{l}{3m}}\frac{1}{v_m^-(K^\circ)}:1\ls m\ls l\Bigr\}.
\end{equation}
Finally, if $K$ has Santaló point at the origin, then
\begin{equation}\label{eq3:prop:analogue-Gian-Milman1}
e_l(K,B_2^n)\ls C\frac{n}{l}\log\left(\frac{en}{l}\right)[1+\|-b(K)\|_K]\sup\Bigl\{2^{-\frac{l}{3m}}\left(\frac{n}{m}\right)^2 w_m(K):1\ls m\ls l\Bigr\}.
\end{equation}
and
\begin{equation*}
e_l(K,B_2^n)\ls C\frac{n}{l}\log\left(\frac{en}{l}\right)[1+\|-b(K)\|_K]\sup\Bigl\{2^{-\frac{l}{3m}}\left(\frac{n}{m}\right)^3\frac{1}{v_m^-(K^\circ)}:1\ls m\ls l\Bigr\}.
\end{equation*}
\end{proposition}
\begin{proof}
Aside from inequality \eqref{eq2:prop:analogue-Gian-Milman1}, we can prove the rest by simply reducing to the symmetric case, covered by Proposition \ref{prop1:Gian-Milman}, and then adjusting for the relevant parameters through the volume comparison results we have seen. Indeed, if $K$ has barycentre at the origin, then we can write
\begin{equation*}
e_l(K,B_2^n)\ls e_l(K-K,B_2^n)\ls C\frac{n}{l}\log\left(\frac{en}{l}\right)\sup\Bigl\{2^{-\frac{l}{3m}}w_m(K-K):1\ls m\ls l\Bigr\}
\end{equation*}
and then recall that
\begin{equation*}
w_m(K-K) \ls C\frac{n}{m}\sup_{x\in {\mathbb R}^n}w_m(K-x)\ls C^\prime\left(\frac{n}{m}\right)^2 w_m(K)
\end{equation*}
by Rudelson's and Fradelizi's results. Finally, the companion inequality to \eqref{eq1:prop:analogue-Gian-Milman1} follows from combining this with Theorem \ref{thm:BS-ineq-proj-sec}.

\smallskip

Similarly, if $K$ has Santaló point at the origin, then we have
\begin{align*}
w_m(K-K)\ls C\frac{n}{m}\sup_{x\in {\mathbb R}^n}w_m(K-x)&\ls C^\prime\left(\frac{n}{m}\right)^2 w_m(K-b(K))
\\
&\ls C^\prime\left(\frac{n}{m}\right)^2 [1+\|-b(K)\|_K]\,w_m(K).
\end{align*}

\bigskip

On the other hand, to prove \eqref{eq2:prop:analogue-Gian-Milman1}, we adapt the proof of Proposition \ref{prop1:Gian-Milman} (essentially repeating it, except that we replace Theorem \ref{thm:Pisier-regular-ellipsoids} by Theorem \ref{main-result}). For the reader's convenience, we sketch the argument. 

Since $l\mapsto e_l(K,B_2^n)$ is decreasing, we can just try to bound $e_l(K,B_2^n)$ for $l$ in some arithmetic progression of integers, and also for $l=1$. 

If $l=1$, observe that $v_1^-(K^\circ)\gs 2r(K^\circ)$ and recall that $\frac{1}{2r(K^\circ)} = \frac{1}{2}R(K)$, while $e_1(K,B_2^n) = R(K)$. Thus the inequality holds for $l=1$ (and in fact it is very crude as $n$ gets larger and larger).

Assume now that $l$ is a multiple of 3. By Theorem \ref{main-result} we can find a $\beta$-regular $M$-ellipsoid ${\cal E}_\beta$ for $K$ for any $\beta\in (0,\frac{2}{5})$ (we will eventually choose $\beta_0=\frac{2}{5}-\frac{1}{5\log(en/l)}$). We can write
\begin{equation*}
e_{l+1}(K,B_2^n)\ls e_{l/3+1}(K,{\cal E}_\beta)e_{2l/3+1}({\cal E}_\beta,B_2^n)\ls D_\beta^{1/\beta}\left(\frac{3n}{\log(2)\,l}\right)^{1/\beta}e_{2l/3+1}({\cal E}_\beta,B_2^n).
\end{equation*}
At the same time, by Lemma \ref{lem:ellipsoid-covering} and Remark \ref{rem:ellipsoid-covering} (i) and (iii), we know that
\begin{equation*}
e_{2l/3+1}({\cal E}_\beta,B_2^n) \ls C\sup_{1\ls m\ls l} 2^{-\frac{2l}{3m}}w_m({\cal E}_\beta) = C\sup_{1\ls m\ls l} 2^{-\frac{2l}{3m}}\frac{1}{v_m^-({\cal E}_\beta^\circ)}\,.
\end{equation*}
But for every $F\in G_{n,m}$ and every $t\gs D_\beta^{1/\beta}$,
\begin{align*}
\frac{{\rm vrad}({\rm Proj}_F(K^\circ))}{{\rm vrad}({\rm Proj}_F({\cal E}_\beta^\circ))} = \frac{|{\rm Proj}_F(K^\circ)|^{1/m}}{|{\rm Proj}_F({\cal E}_\beta^\circ)|^{1/m}}&\ls t\,N({\rm Proj}_F(K^\circ), t{\rm Proj}_F({\cal E}_\beta^\circ))^{1/m}
\\
&\ls t\,N(K^\circ, t{\cal E}_\beta^\circ)^{1/m}\ls t\,\exp(D_\beta n/(m\, t^\beta)),
\end{align*}
and if we choose $t=D_\beta^{1/\beta}\left(\frac{3n}{\log(2)\, l}\right)^{1/\beta}$ we obtain
\begin{equation*}
\frac{{\rm vrad}({\rm Proj}_F(K^\circ))}{{\rm vrad}({\rm Proj}_F({\cal E}_\beta^\circ))}\ls D_\beta^{1/\beta}\left(\frac{3n}{\log(2)\, l}\right)^{1/\beta}2^{l/(3m)}.
\end{equation*}
It follows that
\begin{equation*}
e_{2l/3+1}({\cal E}_\beta,B_2^n)\ls CD_\beta^{1/\beta}\left(\frac{3n}{\log(2)\, l}\right)^{1/\beta}\sup_{1\ls m\ls l} 2^{-\frac{l}{3m}}\frac{1}{v_m^-(K^\circ)}
\end{equation*}
and 
\begin{equation*}
e_{l+1}(K,B_2^n) \ls CD_\beta^{2/\beta}\left(\frac{3n}{\log(2)\, l}\right)^{2/\beta}\sup_{1\ls m\ls l} 2^{-\frac{l}{3m}}\frac{1}{v_m^-(K^\circ)}
\end{equation*}
On setting $\beta= \frac{2}{5}-\frac{1}{5\log(en/l)}$, we complete the proof.
\end{proof}

\begin{proposition}\label{prop:analogue-Gian-Milman2}
Let $K$ be a convex body in ${\mathbb R}^n$. If $K$ has barycentre at the origin, then for every $l\in \{1,2,\ldots,\lfloor\frac{n}{2}\rfloor\}$ we have
\begin{multline}\label{eq1:prop:analogue-Gian-Milman2}
\qquad \quad c_{2l-1}(K,B_2^n) \ls C\left(\frac{n}{l}\right)^3\log\left(\frac{en}{l}\right)\,w_l(K)
\\
\hbox{and} \quad c_{2l-1}(K,B_2^n) \ls C\left(\frac{n}{l}\right)^4\log\left(\frac{en}{l}\right)\frac{1}{v_l^-(K^\circ)}\,. \qquad
\end{multline}
On the other hand, if $K$ has Santaló point at the origin, then, for $l$ in the same range,
\begin{equation}\label{eq2:prop:analogue-Gian-Milman2}
c_{2l-1}(K,B_2^n) \ls C\left(\frac{n}{l}\right)^5\log^2\left(\frac{en}{l}\right)\frac{1}{v_l^-(K^\circ)}\,.
\end{equation}
\end{proposition}
\begin{proof}
Again for \eqref{eq1:prop:analogue-Gian-Milman2} we reduce to the corresponding proposition about the symmetric case, Proposition \ref{prop2:Gian-Milman}:
\begin{equation*}
c_{2l-1}(K,B_2^n) \ls c_{2l-1}(K-K,B_2^n)\ls C\frac{n}{l}\log\left(\frac{en}{l}\right)\,w_l(K-K)\,.
\end{equation*}
As we have seen
\begin{equation*}
w_l(K-K)\ls C\left(\frac{n}{l}\right)^2 w_l(K) \ls C^\prime\left(\frac{n}{l}\right)^3 \frac{1}{v_l^-(K^\circ)}
\end{equation*}
for $K$ centred.

On the other hand, for \eqref{eq2:prop:analogue-Gian-Milman2} we make the necessary adjustments to the proof of Proposition \ref{prop2:Gian-Milman} (by substituting the use of Pisier's Theorem \ref{thm:Pisier-regular-ellipsoids} with Theorem \ref{main-result} and Proposition \ref{prop:Gelfand-numbers}; in fact the proof of that proposition and of \eqref{eq2:prop:analogue-Gian-Milman2} are very similar). 

Fix $l\in \{0, 1,2,\ldots,\lfloor\frac{n-1}{2}\rfloor\}$. By Theorem \ref{main-result} and Proposition \ref{prop:Gelfand-numbers}, for any $\beta \in (0,\frac{2}{5})$ we can find a $\beta$-regular $M$-ellipsoid ${\cal E}_\beta$ for $K$ for which we will have
\begin{equation*}
c_l(K,{\cal E}_\beta)\ls C_0D_\beta^{1/\beta}\left(\frac{n}{l}\right)^{1/\beta}.
\end{equation*}
Again we can assume that $l>1$ since the desired inequality obviously holds for $l=1$ (as $c_1(K,B_2^n) = R(K)=\frac{2}{v_1^-(K^\circ)}$). By the definition of the Gelfand numbers, we can find a subspace $F_0\in G_{n, n-l+1}$ such that
\begin{equation*}
K\cap F_0\subseteq C_0D_\beta^{1/\beta}\left(\frac{n}{l}\right)^{1/\beta}\cdot {\cal E}_\beta\cap F_0\,.
\end{equation*}
Let $\mu_1\gs \mu_2\gs\cdots\gs\mu_{n-l+1}$ be the lengths of the semiaxes of the ellipsoid ${\cal E}_\beta\cap F_0$, and let $F_1\ls F_0$ be the subspace spanned by the $n-2l+1$ shortest of those.
%(thus the orthogonal complement $F_1^\perp \cap F_0$ of $F_1$ within $F_0$ is the subspace spanned by the $l$ longest semiaxes of ${\cal E}_\beta\cap F_0$). 
Then ${\cal E}_\beta\cap F_1\subseteq \mu_{l+1}B_{F_1}$ and
\begin{equation*}
\mu_{l+1} \ls \biggl(\prod_{i=1}^{l}\mu_i\biggr)^{1/l} = {\rm vrad}({\cal E}_\beta\cap F_1^\perp\cap F_0) \ls w_l({\cal E}_\beta) = \frac{1}{v_l^-({\cal E}_\beta^\circ)}.
\end{equation*}
As before, we can compare $v_l^-({\cal E}_\beta^\circ)$ with $v_l^-(K^\circ)$ using the regularity of the covering numbers of $K^\circ$ by ${\cal E}_\beta^\circ$:
\begin{equation*}
\frac{1}{v_l^-({\cal E}_\beta^\circ)}\ls CD_\beta^{1/\beta}\left(\frac{n}{l}\right)^{1/\beta}\,\frac{1}{v_l^-(K^\circ)}.
\end{equation*}
Thus
\begin{equation*}
K\cap F_1 \subseteq C_0D_\beta^{1/\beta}\left(\frac{n}{l}\right)^{1/\beta}\cdot {\cal E}_\beta\cap F_1 \subseteq CD_\beta^{2/\beta}\left(\frac{n}{l}\right)^{2/\beta}\,\frac{1}{v_l^-(K^\circ)}\cdot B_{F_1}.
\end{equation*}
On setting $\beta= \frac{2}{5}-\frac{1}{5\log(en/l)}$, we complete the proof.
\end{proof}

\begin{corollary}{\rm (Corollary to Proposition \ref{prop:analogue-Gian-Milman1})}
Let $K$ be a not-necessarily symmetric {\bf isotropic} convex body in ${\mathbb R}^n$. Then
\begin{equation*}
M(K) \ls C\frac{\log^{\frac{5}{22}}(n)}{n^{\frac{1}{22}}\,L_K}\,.
\end{equation*}
Moreover, we have the `conditional' bound
\begin{equation*}
M(K)\ls C\,[1+h_K(-b(K^\circ))]\,\frac{\log^{\frac{1}{6}}(n)}{n^{\frac{1}{18}}\,L_K}\,.
\end{equation*}
\end{corollary}
\begin{proof}
We combine Dudley's entropy estimate with \eqref{eq2:prop:analogue-Gian-Milman1} from Proposition \ref{prop:analogue-Gian-Milman1} (which we apply with the convex body $K^\circ$, that has Santaló point at the origin), and also the obvious bound $e_l(K^\circ,B_2^n) \ls R(K^\circ) = \frac{1}{r(K)}$:
\begin{equation*}
\sqrt{n}M(K) = \sqrt{n}M^\ast(K^\circ) \ls C\sum_{l=1}^n\frac{1}{\sqrt{l}}\min\left\{\frac{1}{r(K)},\,\left(\frac{n}{l}\right)^5\log^2\left(\frac{en}{l}\right)\sup_{1\ls m\ls l}\Bigl\{2^{-\frac{l}{3m}}\frac{1}{v_m^-(K)}\Bigr\}\right\}.
\end{equation*}
By the properties of the isotropic position that we listed in Subsection \ref{subsec:properties-isotropic}, we have $r(K)\gs c_0L_K$, and moreover, for every $F\in G_{n,m}$,
%\begin{equation}\label{eq:bounds4proj}
%{\rm vrad}({\rm Proj}_F(K)) \gs \frac{1}{|B_2^m|^{1/m}}\frac{1}{|K\cap F^\perp|^{1/m}} \gtrsim \sqrt{m}\,\frac{L_K}{\sup_{M\subset {\mathbb R}^m} L_M}
%\end{equation}
\begin{equation*}
{\rm vrad}({\rm Proj}_F(K)) \gtrsim \sqrt{m}\,\frac{L_K}{\sup_{M\subset {\mathbb R}^m} L_M}\,.
\end{equation*}
%by standard properties of the isotropic position (that we recalled in Section \ref{sec:prelim}).
By using the latest progress on the bounds for $\sup_{M\subset {\mathbb R}^m} L_M$, we obtain
\begin{equation*}
\frac{1}{v_m^-(K)} \ls C^\prime \frac{1}{\sqrt{m}} \frac{\sqrt{\log(m)}}{L_K}\ls C^\prime \frac{1}{\sqrt{m}} \frac{\sqrt{\log(n)}}{L_K},
\end{equation*}
and thus we can compute
\begin{equation*}
\sup_{1\ls m\ls l}\Bigl\{2^{-\frac{l}{3m}}\frac{1}{v_m^-(K)}\Bigr\} \ls C^{\prime\prime} \frac{1}{\sqrt{l}} \frac{\sqrt{\log(n)}}{L_K}.
\end{equation*}
By bounding $\log^2(en/l)$ by $\log^2(n)$ too for simplicity, we can finally conclude that
\begin{align*}
\sqrt{n}M(K) \ls C\sum_{l=1}^n\min\left\{\frac{1}{\sqrt{l}\,L_K},\,\frac{n^5}{l^6}\frac{\log^{5/2}(n)}{L_K}\right\}
\ls C^\prime\, \frac{n^{\frac{5}{11}}\,\log^{\frac{5}{22}}(n)}{L_K}.
\end{align*}

Similarly we deduce the other bound, by using the inequality after \eqref{eq3:prop:analogue-Gian-Milman1} this time.
\end{proof}

\bigskip

\noindent {\bf Acknowledgement.} The author would like to thank Alexander Litvak for suggesting the use of an `average' ellipsoid in the proof of Theorem \ref{main-result}.

\medskip
\bigskip

{\footnotesize

\noindent \textsc{Beatrice-Helen Vritsiou:} Department of Mathematical and Statistical Sciences,
University of Alberta, CAB 632, Edmonton, AB, Canada T6G 2G1

\smallskip

\noindent 
{\it E-mail:} \texttt{vritsiou@ualberta.ca}

}

\end{document}